\documentclass[]{interact}

\usepackage{mathtools}
\usepackage{bm}
\usepackage{array}
\graphicspath{{figures/}}
\usepackage{placeins}
\usepackage{enumitem}
\setlist[enumerate]{leftmargin=.5in}
\setlist[itemize]{leftmargin=.5in}
\usepackage{tikz}
\usetikzlibrary{arrows.meta,calc}
\usepackage{microtype}

\usepackage[numbers,sort&compress]{natbib}
\bibpunct[, ]{[}{]}{,}{n}{,}{,}

\makeatletter
\def\NAT@def@citea{\def\@citea{\NAT@separator}}
\makeatother

\usepackage[colorlinks=true,allcolors=blue]{hyperref}

\numberwithin{equation}{section}

\theoremstyle{plain}
\newtheorem{theorem}{Theorem}[section]
\newtheorem{lemma}[theorem]{Lemma}

\newtheorem{proposition}[theorem]{Proposition}
\theoremstyle{definition}
\newtheorem{definition}[theorem]{Definition}
\theoremstyle{remark}
\newtheorem{remark}[theorem]{Remark}
\newtheorem{assumption}[theorem]{Assumption}

\DeclareMathOperator{\ran}{Ran}
\DeclareMathOperator{\dist}{dist}
\DeclareMathOperator{\Var}{Var}

\DeclareMathOperator{\Real}{Re}
\DeclareMathOperator{\Imag}{Im}
\DeclareMathOperator{\sgn}{sgn}

\newcommand{\R}{\mathbb{R}}
\newcommand{\C}{\mathbb{C}}

\newcommand{\Sone}{\mathbb{S}^{1}}
\newcommand{\LL}{L^2(\Sone)}
\newcommand{\E}{\mathbb{E}}
\newcommand{\Prob}{\mathbb{P}}
\newcommand{\norm}[1]{\left\lVert #1\right\rVert}
\newcommand{\opnorm}[1]{\left\lVert #1\right\rVert_{\mathrm{op}}}
\newcommand{\ip}[2]{\left\langle #1,#2\right\rangle}
\newcommand{\Fh}{F_{\#}}

\newcommand{\cS}{\mathcal S}
\newcommand{\cQ}{\mathcal Q}

\newcommand{\dd}{\,\mathrm d}

\definecolor{seablue}{RGB}{38,70,130}
\definecolor{skyfill}{RGB}{233,241,251}
\definecolor{obsred}{RGB}{200,52,52}
\definecolor{obsfill}{RGB}{250,206,196}
\definecolor{senteal}{RGB}{17,122,118}
\definecolor{wavecol}{RGB}{233,140,38}
\definecolor{graypr}{RGB}{150,150,160}

\begin{document}

\title{Pathwise Factorization and Uncertainty Quantification for
Impedance Obstacles with Random Boundaries}

\author{
\name{George Pelekanos\textsuperscript{a}\thanks{CONTACT George Pelekanos.
Email: gpeleka@siue.edu} and Koung Hee Leem\textsuperscript{a}}
\affil{\textsuperscript{a}Department of Mathematics and Statistics,
Southern Illinois University Edwardsville, Edwardsville, IL 62026, USA}
}

\maketitle

\begin{abstract}
We develop a factorization method with finite-sample statistical
guarantees for two-dimensional acoustic scattering by impedance obstacles
whose boundary is random and whose data are finitely many noisy multistatic
measurements.  In the two-dimensional far-field normalization, the
symmetrization must act on a phase-rotated far-field operator; we prove the
optical and range identities for every realization of a compact family of
uniformly bi-Lipschitz domains, with jointly measurable regularized
indicators.  A quadrature-weighted finite-sample analysis separates angular
approximation, measurement noise, test-vector error, and regularization,
and yields a consistency theorem for empirical means and variances, joint
in all data parameters.  Bounded indicators admit simultaneous confidence
bands, a diagonal threshold recovers the probability that a sampling point
lies in the random obstacle, and a nested design yields unbiased variance
components.  Computations use a combined-field integral equation proved
invertible at all frequencies, and the finite-grid guarantees are examined
in a coverage study against a high-accuracy Monte Carlo reference.
\end{abstract}


\begin{amscode}
35R30, 35R60, 47A55, 65N21, 60H35, 62G20
\end{amscode}

\section{Introduction}\label{sec:introduction}

Qualitative inverse-scattering methods recover an obstacle without solving a
nonlinear optimization problem and without a close initial guess.  Their
classical theory nevertheless assumes an exact continuum data operator and a
fixed obstacle.  In practice, the boundary may vary because of manufacturing
tolerances, deformation, erosion, biological variability, or registration
error, while the data are available only at finitely many directions and are
contaminated by measurement noise.  This paper treats these effects together
for an absorbing impedance obstacle.

The boundary condition
\begin{equation}\label{eq:intro-robin}
 \partial_\nu u+ik\lambda u=0,
 \qquad \lambda\ge \lambda_0>0,
\end{equation}
models a passive coating or acoustic lining.  With the convention used here,
$\lambda$ is real and positive; this is already an absorbing law.  The
far-field operator need not be normal for a nonsymmetric impedance obstacle,
so a normal-operator range test cannot be assumed.  Nonnormality is not forced,
however; the impedance disc is diagonal in the Fourier basis and is therefore
normal.  The general construction is the Kirsch--Grinberg
symmetrization together with the generalized boundary-condition factorization
of Chamaillard, Chaulet, and Haddar \cite{Gri06,KG08,CCH13}.  In two
dimensions there is an additional point of normalization.  The far-field
convention used here carries the Hankel phase $e^{-i\pi/4}$, and the
symmetrization must be applied to the phase-rotated far-field operator.
Section~\ref{sec:factorization} proves the exact optical identity behind
this normalization, and gives an obstacle in the admissible class for which
the unrotated symmetrization fails.

Two further structural points require care.  First, the factorization is
built on the data-to-pattern operator, which maps impedance boundary data
to far-field patterns, rather than on the full Herglotz Cauchy trace.  Its
Cauchy pairs range only over Cauchy data of interior
Helmholtz solutions, which are not dense in the product of the natural
trace spaces; a product-space factorization relying on dense range is
therefore invalid.  Moreover, in the physical far-field amplitude used here, the
middle operator of the symmetric factorization carries a positive real
scalar and the two-dimensional Hankel phase; both are immaterial for ranges
but essential in the proofs, and Section~\ref{sec:factorization} transfers
the standard normalization explicitly and obtains the range identity of
the factorization method.  Second, the continuum far-field operator acts on
square-integrable functions of the observation direction, while a raw
multistatic matrix acts on Euclidean vectors; the two are related by
angular quadrature weights, whose omission makes the matrix norm and the
regularization parameter depend on the direction count.  All theory and
computations therefore use quadrature-weighted data matrices and test
vectors; for uniform full-circle sampling with $M$ directions, the weighted
regularization parameter is $2\pi/M$ times the raw one.

This work sits at the junction of three strands of earlier work.  The
factorization method originates with Kirsch \cite{Kir98,Kir99} and is
developed in the monograph of Kirsch and Grinberg \cite{KG08}; its extension
to impedance and to general boundary conditions, on which
Section~\ref{sec:factorization} relies, is due to Grinberg \cite{Gri06} and
to Chamaillard, Chaulet, and Haddar \cite{CCH13}.  Regularization of the
method under deterministic, worst-case data error was analyzed by Lechleiter
\cite{Lec06} and, for perturbed positive compact operators, by Harris
\cite{Har23}; neither treats a statistical noise model or finite angular
sampling.  Randomness has entered qualitative inversion along other axes.
Griesmaier and Raumer \cite{GR22} factorize the data covariance operator for
random \emph{sources}, assuming the exact covariance is available, and
Garnier, Haddar, and Montanelli \cite{GHM23} develop a linear sampling
method for randomly distributed point sources from cross-correlated passive
measurements.  Barth,
Harrach, Hyv\"onen, and Mustonen \cite{BHHM17} show in electrical impedance
tomography that the method applied to the \emph{mean} data operator detects a
fixed inclusion with random \emph{material}; averaging the operator is
legitimate there precisely because the geometry is deterministic.  For random
\emph{geometry}, the domain-mapping approach is well developed on the forward
side, in shape uncertainty quantification for Helmholtz scattering
\cite{HSSS18} and for elliptic problems on random domains \cite{HPS16}, and
recent work on the inverse side reconstructs random obstacles by iterative
multi-frequency optimization \cite{SXL26}.  A qualitative method for random
boundaries --- with a range characterization holding for each realization and
statistical guarantees for finitely many noisy measurements --- appears not
to have been available, and supplying it is the purpose of this paper.

The contributions are the following.  On the deterministic side, we track the
standard and physical two-dimensional far-field normalizations and obtain the
correct rotated factorization with a positive real middle-operator scaling.
We prove that neither the direct nor the adjoint interior impedance problem has
real eigenvalues, transfer the exact optical identity with its explicit
constant, and derive the pathwise range theorem through the data-to-pattern
operator.  The random geometry is placed on a compact parametric family of
uniformly bi-Lipschitz domains.  A fixed-domain variational pullback yields
operator-norm continuity of the forward map, joint measurability of the
regularized indicator, and an extended-real measurable exact Picard field.

On the finite-sample side, a quadrature-weighted perturbation theorem
separates angular approximation, the error from repeated noisy acquisitions
of the multistatic matrix (measurement frames), test-vector error,
and regularization.  Deterministic confidence schedules give pathwise
convergence in probability or almost surely.  A uniform random-geometry
version yields a consistency theorem for empirical means and variances,
joint in the direction count, the frame count, the number of geometry
draws, and the regularization parameter.  The regularization rates that
these estimates call for, together with the source conditions and sharp
spectral asymptotics that determine them, are studied separately
in \cite{PartII}.

For uncertainty quantification, the bounded reciprocal indicator satisfies
simultaneous finite-grid concentration inequalities, upgraded by a spatial
Lipschitz estimate to continuum confidence bands for the mean regularized
indicator; the bands are not, by themselves, confidence regions for
obstacle membership.
A separate diagonal threshold construction converges to the probability
that a sampling point lies in the random obstacle, together with its
Bernoulli variance.  A nested
geometry--measurement design --- repeated reconstructions grouped within
each geometry draw --- gives unbiased, strongly consistent variance
components with finite-sample bounds and an acquisition-cost rule for the
number of repeats per geometry.

On the computational side, physical well-posedness is separated from
representation-induced resonances.  A combined-field family is proved
invertible at every frequency; its discretization is validated against the
analytic disc and generates all reported ensembles, with the conditioning of
every realization recorded.  Ensemble-wide
condition-number diagnostics on the same random paths reveal several poorly
conditioned single-layer solves, whereas the combined-field systems are
well conditioned on every sampled path; the single-layer representation is retained for
validation and for that comparison.

The uncertainty quantified here is the push-forward of a prescribed geometry
law through the pathwise reconstruction map.  It is not a Bayesian posterior
conditioned on one observed data set \cite{Stuart10}.  The distinction is
revisited in Section~\ref{sec:conclusion}.

The paper is organized as follows.  Section~\ref{sec:problem} formulates
the scattering problem and the admissible random-domain class;
Section~\ref{sec:factorization} establishes the normalization-correct
factorization, the optical identity, the range theorem, and the regularized
Picard indicator; Section~\ref{sec:pathwise} proves pathwise validity and
measurability; Section~\ref{sec:discrete} develops the quadrature-weighted
finite-sample theory; Section~\ref{sec:uq} gives the confidence bands, the
joint consistency theorem, the nested variance decomposition, and
occupancy recovery; Section~\ref{sec:numerics} reports the computations;
Section~\ref{sec:conclusion} concludes.  Validation of the discretizations and the complete experiment
specification are collected in Appendices~\ref{app:validation}
and~\ref{app:params}, and the
regularization rates with the associated sharp disc asymptotics are
treated separately in \cite{PartII}.  The logical
spine is short.  The deterministic range identity of
Section~\ref{sec:factorization} feeds the pathwise measurability of
Section~\ref{sec:pathwise}; the discrete perturbation theory of
Section~\ref{sec:discrete} turns it into finite-sample statements; and the
Monte Carlo analysis of Section~\ref{sec:uq} rests on these two
 ingredients alone.  Sections~\ref{sec:factorization}
and~\ref{sec:pathwise} establish the deterministic and pathwise results for
every member of the admissible class.  The finite-data results of
Section~\ref{sec:discrete} are conditional on the stated quadrature,
test-vector, norm, and measurement-noise hypotheses, and the joint
consistency theorem additionally assumes that these hold uniformly over the
geometry family; the coverage statements of Section~\ref{sec:numerics} are
numerical observations.

\section{Scattering problem and admissible random domains}\label{sec:problem}

\subsection{Fixed obstacle and far-field operator}
Let $D\subset\R^2$ be a bounded open scatterer with finitely many
connected components, $C^2$ boundary $\Gamma=\partial D$, and connected
exterior.  The $C^2$ assumption is used for
the parameter-continuity and numerical sections; the deterministic
factorization itself extends to the Lipschitz setting.  The unit normal $\nu$
points into the exterior.  For $d\in\Sone$, the incident wave is
$u^i(x,d)=e^{ikx\cdot d}$, and the total field $u=u^i+u^s$ satisfies
\begin{equation}\label{eq:forward}
 \Delta u+k^2u=0\quad\text{in }\R^2\setminus\overline D,
 \qquad
 \partial_\nu u+ik\lambda u=0\quad\text{on }\Gamma,
\end{equation}
with $k>0$ and a real impedance
$\lambda\in L^\infty(\Gamma)$ satisfying
$0<\lambda_0\le\lambda\le\lambda_1$.  The scattered field is outgoing and has
the asymptotic form \cite{CK19}
\begin{equation}\label{eq:far-asymptotic}
 u^s(x,d)=\frac{e^{ik|x|}}{|x|^{1/2}}
 \left(u^\infty(\widehat x,d)+O(|x|^{-1})\right),
 \qquad \widehat x=x/|x|.
\end{equation}
The far-field operator $F:\LL\to\LL$ is
\begin{equation}\label{eq:Fdef}
 (Fg)(\widehat x)=\int_{\Sone}u^\infty(\widehat x,d)g(d)\,\dd s(d).
\end{equation}
All complex Hilbert-space inner products are linear in the first argument.

\begin{figure}[tbhp]
\centering
\begin{tikzpicture}[scale=0.88,>=Stealth]
\fill[skyfill,rounded corners=2pt] (-2.1,-2.1) rectangle (2.1,2.1);
\draw[graypr,rounded corners=2pt] (-2.1,-2.1) rectangle (2.1,2.1);
\node[graypr,font=\footnotesize] at (1.82,-1.87) {$K$};
\fill[obsfill,opacity=0.85] plot[domain=0:360,samples=140,smooth cycle]
  ({\x}:{0.62*(1.45+0.34*cos(3*\x))});
\draw[obsred,very thick] plot[domain=0:360,samples=140,smooth cycle]
  ({\x}:{0.62*(1.45+0.34*cos(3*\x))});
\draw[senteal,dashed,thick] plot[domain=0:360,samples=140,smooth cycle]
  ({\x}:{0.62*(1.45+0.34*cos(3*\x)+0.19*cos(\x)-0.13*sin(2*\x))});
\draw[graypr,dashed,thick] plot[domain=0:360,samples=140,smooth cycle]
  ({\x}:{0.62*(1.45+0.34*cos(3*\x)-0.16*cos(2*\x)+0.12*sin(\x))});
\node[obsred!70!black,font=\footnotesize] at (0,0.1) {$D_{Y}$};
\node[obsred!70!black,font=\footnotesize] at (-0.15,-1.62)
  {$\partial_\nu u+ik\lambda u=0$};
\draw[graypr,thin] (-0.32,-1.45) -- (-0.38,-1.04);
\draw[seablue,dashed] (0,0) circle (2.55);
\foreach \a in {0,30,...,330} \fill[seablue] (\a:2.55) circle (0.035);
\draw[wavecol,very thick,->] (-3.65,1.95) -- (-2.1,1.1);
\draw[wavecol,thick] (-3.28,2.29) -- (-3.62,1.67);
\draw[wavecol,thick] (-2.98,2.12) -- (-3.31,1.51);
\draw[wavecol,thick] (-2.67,1.95) -- (-3.00,1.34);
\node[wavecol!80!black,font=\footnotesize] at (-3.15,2.66) {$u^{i}=e^{ikx\cdot d}$};
\node[wavecol!80!black,font=\footnotesize] at (-2.42,0.82) {$d$};
\foreach \a in {35,-20,215} \draw[seablue,->] (\a:1.62) -- (\a:2.45);
\node[seablue,font=\footnotesize] at (35:3.15) {$u^{\infty}(\widehat x,d)$};
\node[seablue,font=\footnotesize] at (-20:2.92) {$\widehat x$};
\node[seablue,font=\footnotesize] at (0,-2.88) {$M$ directions};
\fill[black] (-1.4,1.4) circle (0.045);
\node[font=\footnotesize] at (-1.38,1.72) {$z,\ \varphi_z$};
\end{tikzpicture}
\caption{The setting.  One realization $D_Y$ of the random impedance
obstacle (solid, with the boundary condition it carries) and two further
realizations of the random boundary (dashed).  Plane waves are sent from,
and far-field patterns observed in, $M$ full-aperture directions on
$\Sone$ (dotted circle); the indicator is evaluated at sampling points $z$
of the window $K$ through the test function $\varphi_z$.}
\label{fig:setup}
\end{figure}
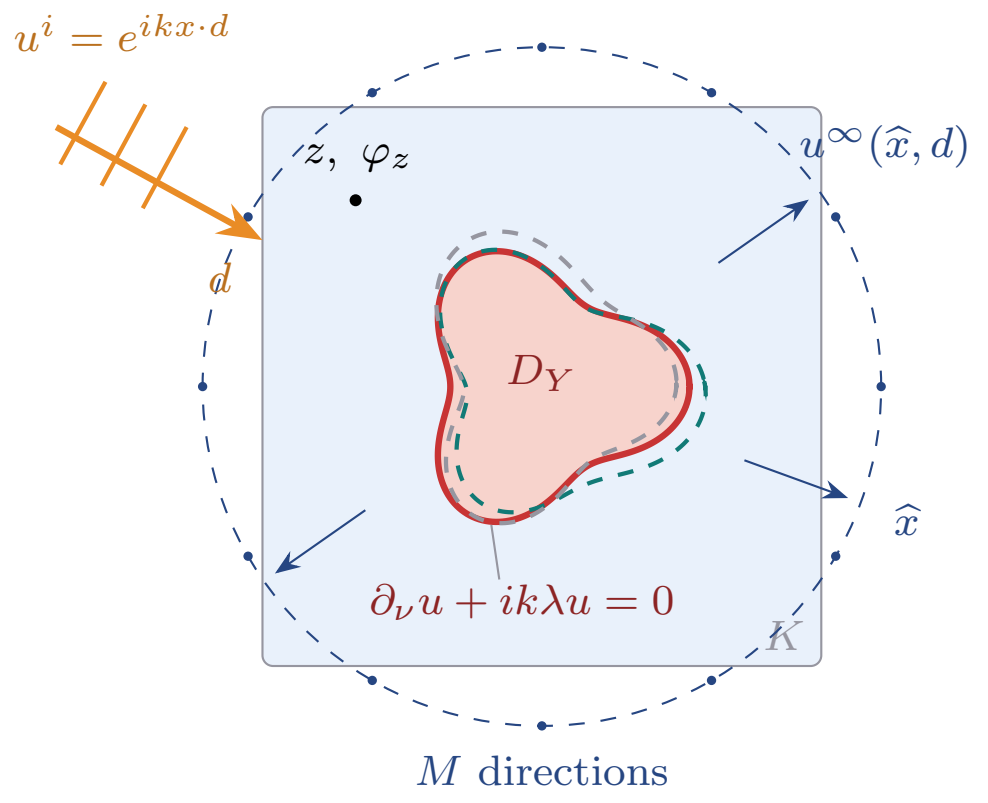

For $g\in\LL$, define the Herglotz wave
\begin{equation}\label{eq:herglotz}
 v_g(x)=\int_{\Sone}e^{ikx\cdot d}g(d)\,\dd s(d),
 \qquad H_0g=\gamma_0v_g,
 \qquad H_1g=\gamma_1v_g.
\end{equation}
The sampling function is
\begin{equation}\label{eq:test-function}
 \varphi_z(\widehat x)=e^{-ik\widehat x\cdot z},
 \qquad z\in\R^2.
\end{equation}
A nonzero normalization constant in \eqref{eq:test-function} does not affect any
range test.  Figure~\ref{fig:setup} summarizes the setting.

\subsection{Compact parametric random geometry}
Let $\widehat D\subset\R^2$ be a fixed bounded reference obstacle with
finitely many components, $C^2$ boundary, connected exterior, and positive
separation between distinct components, and let $Y$ be a compact metric
parameter space.

\begin{assumption}[Admissible parametric family]\label{ass:parametric}
For every $y\in Y$, let $\Psi_y:\R^2\to\R^2$ be a global $C^2$
diffeomorphism with the following properties.
\begin{enumerate}[label=(A\arabic*),leftmargin=2.7em]
\item $y\mapsto\Psi_y$ is continuous in $C^2$ on a fixed ball containing the
supports of all $\Psi_y-\mathrm{Id}$;
\item $\Psi_y$ is the identity outside a larger fixed ball;
\item there are constants $0<c_\Psi\le C_\Psi<\infty$ with
\[
 c_\Psi|x-x'|\le |\Psi_y(x)-\Psi_y(x')|
 \le C_\Psi|x-x'|
\]
for all $x,x'$ and all $y$;
\item $D_y=\Psi_y(\widehat D)$ has connected exterior, and distinct reference
components have a uniform positive separation after mapping;
\item the impedance is deterministic and transported with the material, so that
$\lambda_y\circ\Psi_y=\widehat\lambda$, where
$\widehat\lambda\in L^\infty(\partial\widehat D)$ satisfies
$0<\lambda_0\le\widehat\lambda\le\lambda_1$.
\end{enumerate}
Let $Y(\omega)$ be a measurable $Y$-valued random variable.
\end{assumption}

The following elementary construction is useful both theoretically and
numerically.

\begin{lemma}[Near-identity domain maps]\label{lem:diffeomorphism}
Let $V_y\in C_c^2(\R^2;\R^2)$ depend continuously on $y\in Y$ and suppose
\[
 \sup_{y\in Y}\norm{DV_y}_{L^\infty}\le q<1.
\]
Then $\Psi_y=\mathrm{Id}+V_y$ is a global bi-Lipschitz $C^2$
diffeomorphism and
\[
 (1-q)|x-x'|\le|\Psi_y(x)-\Psi_y(x')|
 \le(1+q)|x-x'|.
\]
In particular, simplicity of each boundary component and positive component
separation are preserved uniformly.
\end{lemma}

\begin{proof}
The displayed bounds follow from the mean-value theorem.  They give
injectivity and properness because $V_y$ is compactly supported.  The Jacobian
$I+DV_y$ is invertible everywhere by the Neumann series.  The global inverse
function theorem then gives a global $C^2$ diffeomorphism.  The lower Lipschitz
bound scales every inter-component distance by at least $1-q$.
\end{proof}

\begin{remark}[Compactly supported extension of an affine prior]\label{rem:affine-extension}
Suppose an affine map $x\mapsto A_yx+b_y$ is used on a ball $B_R$ containing
the obstacle, with $\sup_y\norm{A_y-I}<1/3$, $\sup_y|b_y|<\infty$, and
$y\mapsto(A_y,b_y)$ continuous on the compact set $Y$.  Put $a=\sup_y\norm{A_y-I}<1/3$ and $b=\sup_y|b_y|<\infty$.  Choose
$\zeta\in C^\infty(\R)$ equal to one on $(-\infty,0]$ and zero on
$[1,\infty)$ with $\norm{\zeta'}_\infty\le2$, put
$\chi_L(x)=\zeta\bigl((|x|-R)/L\bigr)$, and set
$V_y(x)=\chi_L(x)((A_y-I)x+b_y)$.  Then
\[
 \sup_y\norm{DV_y}_\infty
 \le a+\frac{2}{L}\bigl(a(R+L)+b\bigr)
 \longrightarrow 3a<1
 \qquad(L\to\infty),
\]
so $\sup_y\norm{DV_y}_\infty<1$ for all sufficiently large $L$; the
threshold $1/3$ in the hypothesis is exactly what makes the limit
admissible.  Lemma~\ref{lem:diffeomorphism} therefore
produces a global compactly supported diffeomorphism that agrees with the
affine map on every boundary used in the computation.
\end{remark}

\section{The deterministic impedance factorization}\label{sec:factorization}

\subsection{Far-field normalization and the data-to-pattern factorization}
\label{sec:normalization-factorization}
Set
\[
 Z_k=ikM_\lambda:H^{1/2}(\Gamma)\longrightarrow H^{-1/2}(\Gamma),
\]
where $M_\lambda$ denotes multiplication by $\lambda$.  For merely
$\lambda\in L^\infty(\Gamma)$ this map is bounded through
$H^{1/2}(\Gamma)\hookrightarrow L^2(\Gamma)\xrightarrow{M_\lambda}
L^2(\Gamma)\hookrightarrow H^{-1/2}(\Gamma)$; no stronger multiplier
regularity is used.  For
$f\in H^{-1/2}(\Gamma)$, let $w_f$ be the unique outgoing solution of
\begin{equation}\label{eq:Gproblem}
 \Delta w_f+k^2w_f=0\quad\text{in }\R^2\setminus\overline D,
 \qquad
 \partial_\nu w_f+Z_kw_f=f\quad\text{on }\Gamma,
\end{equation}
and define the data-to-pattern operator in the physical far-field
normalization \eqref{eq:far-asymptotic} by
\begin{equation}\label{eq:Gdef}
 G:H^{-1/2}(\Gamma)\longrightarrow\LL,
 \qquad Gf=w_f^\infty.
\end{equation}
The operator $G$ is compact and injective and has dense range
\cite[Chapter~2]{KG08}.  Put
\[
 B:=Z_kH_0+H_1.
\]
The scattered field generated by the Herglotz incident wave $v_g$ has boundary
datum $-Bg$, and hence
\begin{equation}\label{eq:first-factor}
 F=-GB.
\end{equation}

The normalization must be tracked in the second factorization.  The generalized
factorization theorem of \cite{CCH13} uses the standard two-dimensional
far-field convention
\[
 u^s(x)=\gamma_2\frac{e^{ik|x|}}{|x|^{1/2}}
 \bigl(u^\infty_{\rm st}(\widehat x)+O(|x|^{-1})\bigr),
 \qquad
 \gamma_2=\frac{e^{i\pi/4}}{\sqrt{8\pi k}}.
\]
Let $F_{\rm st}$ and $G_{\rm st}$ denote the corresponding far-field and
data-to-pattern operators.  Comparing this convention with
\eqref{eq:far-asymptotic} gives
\begin{equation}\label{eq:normalization-transfer}
 F=\gamma_2F_{\rm st},
 \qquad
 G=\gamma_2G_{\rm st}.
\end{equation}
For the middle operator $T:H^{1/2}(\Gamma)\to H^{-1/2}(\Gamma)$ used in
\cite{CCH13}, the standard-normalized adjoint identity is $B^*=G_{\rm st}T$.
Consequently,
\begin{equation}\label{eq:adjoint-factor}
 B^*=\gamma_2^{-1}GT,
 \qquad
 F=-\overline{\gamma_2}^{\,-1}GT^*G^*.
\end{equation}
For the local impedance considered here, $T$ is the specialization of
\[
 T=Z_kS_kZ_k^*+N_k+Z_kK_k+K_k'Z_k^*,
\]
where $S_k,K_k,K_k'$, and $N_k$ are the single-layer, double-layer,
adjoint-double-layer, and hypersingular boundary operators, respectively;
the same symbols denote the corresponding layer potentials when applied to a
density and evaluated off $\Gamma$.
The mapping, Fredholm, and sign properties used below are established in
\cite{CCH13}.

Since $\overline{\gamma_2}^{\,-1}=\sqrt{8\pi k}\,e^{i\pi/4}$, the rotated
operator
\begin{equation}\label{eq:rotated}
 \widetilde F=e^{-i\pi/4}F
\end{equation}
admits the factorization
\begin{equation}\label{eq:rotated-factor}
 \widetilde F
 =-\sqrt{8\pi k}\,GT^*G^*
 =-G\widetilde T^{\,*}G^*,
 \qquad
 \widetilde T=\sqrt{8\pi k}\,T.
\end{equation}
The multiplier in $\widetilde T$ is positive and real, so it preserves the
coercivity-up-to-compactness and Fredholm structure of $T$.  In particular,
no complex rotation of the middle operator is required.

\begin{remark}[Why the full Cauchy trace is not the factorization space]
\label{rem:not-dense}
The map
\[
 g\longmapsto(\gamma_0v_g,\gamma_1v_g)
 \in H^{1/2}(\Gamma)\times H^{-1/2}(\Gamma)
\]
does not have dense range in the product space.  In this remark,
$\langle\cdot,\cdot\rangle_{H^{-1/2},H^{1/2}}$ denotes the \emph{bilinear}
duality pairing extending the boundary integral
$\int_\Gamma\phi\,\psi\,\dd s$, without complex conjugation, so that Green's
second identity applies verbatim.  Let $w$ be any nonzero interior Helmholtz
solution.  The continuous linear functional
\[
 \ell(f,q)
 =\bigl\langle q,\gamma_0w\bigr\rangle_{H^{-1/2},H^{1/2}}
 -\bigl\langle\gamma_1w,f\bigr\rangle_{H^{-1/2},H^{1/2}},
 \qquad
 (f,q)\in H^{1/2}(\Gamma)\times H^{-1/2}(\Gamma),
\]
vanishes on every Herglotz Cauchy pair by Green's second identity.  Moreover
$\ell\ne0$.  If $\gamma_0w$ and $\gamma_1w$ both vanished, $w$ would have zero
Cauchy data and hence $w\equiv0$ by Cauchy uniqueness.  The closure of the
Cauchy-trace range is therefore contained in the proper closed subspace
$\ker\ell$.  This obstruction is precisely why the factorization is built on
$G$ rather than on an unconstrained product trace space.
\end{remark}

\subsection{Direct and adjoint interior impedance problems}
Since $\lambda$ is real,
\begin{equation}\label{eq:Zadjoint}
 Z_k^*=-ikM_\lambda.
\end{equation}
In \cite{CCH13}, the interior non-eigenvalue hypothesis of the generalized
factorization theorem is stated for the impedance operator $Z_k$
(Definition~4.8 and Theorem~4.11 of the arXiv version); its proof enters
through injectivity of the adjoint factor $B^*$, for which the natural
interior problem carries $Z_k^*$.  For the scalar real impedance considered
here the two interior problems are complex conjugates of one another, so the
two formulations coincide, and positive real impedance removes both possible
exceptional sets at once.

\begin{lemma}[No real direct or adjoint impedance eigenvalues]\label{lem:noeig}
Let $\lambda$ be as in \eqref{eq:forward}, so that $\lambda$ is real with
$0<\lambda_0\le\lambda\le\lambda_1$, and let $u\in H^1(D)$ solve
\begin{equation}\label{eq:interior-both}
 \Delta u+k^2u=0\quad\text{in }D,
 \qquad
 \partial_\nu u\pm ik\lambda u=0\quad\text{on }\Gamma.
\end{equation}
Then $u=0$.  Thus $k^2$ is neither an interior eigenvalue for $Z_k$ nor for
$Z_k^*$.
\end{lemma}

\begin{proof}
Green's identity gives
\[
 \int_D\bigl(|\nabla u|^2-k^2|u|^2\bigr)\,\dd x
 =\int_\Gamma\overline u\,\partial_\nu u\,\dd s
 =\mp ik\int_\Gamma\lambda|u|^2\,\dd s.
\]
The left side is real and the right side is purely imaginary, so both
vanish.  Since $\lambda$ is of one sign with $\lambda\ge\lambda_0>0$, the
identity $\int_\Gamma\lambda|u|^2\dd s=0$ forces $u|_\Gamma=0$; no
cancellation is available, which is where positivity rather than mere
nonvanishing of $\lambda$ is used.  The boundary condition then
implies $\partial_\nu u|_\Gamma=0$.  Cauchy uniqueness for the Helmholtz
equation yields $u=0$.
\end{proof}

\subsection{Optical identity and exact symmetrization}
For a bounded operator $A$ on a complex Hilbert space, write
\[
 \Real A=\frac12(A+A^*),
 \qquad
 \Imag A=\frac{1}{2i}(A-A^*),
 \qquad
 |H|=(H^2)^{1/2}\quad(H=H^*).
\]
The phase in \eqref{eq:rotated} is forced by the two-dimensional Hankel
asymptotics.  The following proof transfers the standard operator optical
identity and therefore avoids pointwise stationary-phase formulas for an
arbitrary $L^2$ Herglotz density.

\begin{lemma}[Exact optical identity and strict positivity]\label{lem:optical}
For the far-field normalization \eqref{eq:far-asymptotic} and every
$g\in\LL$,
\begin{equation}\label{eq:optical}
 \ip{\Imag \widetilde Fg}{g}_{\LL}
 =\sqrt{\frac{k}{8\pi}}
 \left(\norm{Fg}_{\LL}^2
 +\int_\Gamma\lambda\,|u_g|^2\,\dd s\right),
\end{equation}
where $u_g=v_g+u_g^s$ is the total field for Herglotz incidence with kernel
$g$.  Consequently $\Imag\widetilde F\succeq0$, and
\begin{equation}\label{eq:strict-ImF}
 \ip{\Imag\widetilde Fg}{g}_{\LL}>0\qquad(g\ne0).
\end{equation}
\end{lemma}

\begin{proof}
In the standard normalization, the generalized optical identity
\cite[Lemma~4.10 and equation~(19) of the arXiv version]{CCH13},
specialized to
$Z_k=ikM_\lambda$, is
\[
 \ip{\Imag F_{\rm st}g}{g}
 =k|\gamma_2|^2\norm{F_{\rm st}g}^2
 +k\int_\Gamma\lambda|u_g|^2\,\dd s.
\]
By \eqref{eq:normalization-transfer},
$F_{\rm st}=\sqrt{8\pi k}\,\widetilde F$ and
$\norm{\widetilde Fg}=\norm{Fg}$.  Dividing the displayed identity by
$\sqrt{8\pi k}$ yields \eqref{eq:optical}.

Suppose the right side of \eqref{eq:optical} vanishes.  Then $Fg=0$, so the
scattered Herglotz field has vanishing far field and is zero in the exterior
by Rellich's lemma.  The absorption term and $\lambda\ge\lambda_0$ give
$u_g|_\Gamma=0$, and the boundary condition gives
$\partial_\nu u_g|_\Gamma=0$.  The incident Herglotz wave therefore has zero
Cauchy data on $\Gamma$ and, being entire, vanishes identically by unique
continuation.  Injectivity of the Herglotz map gives $g=0$, proving
\eqref{eq:strict-ImF}.
\end{proof}

\begin{remark}[The rotation cannot be dropped]\label{rem:rotation}
$\Imag F$ itself can be indefinite within the present impedance class, so
the phase in \eqref{eq:rotated} does real work.  For the impedance disc of
radius $a$ the eigenvalues of $F$ are
$\kappa_m=\sqrt{8\pi/k}\,e^{-i\pi/4}a_m$, where
\begin{equation}\label{eq:disc-am}
 a_m=-\frac{J_m'(ka)+i\lambda J_m(ka)}
 {H_m^{(1)\prime}(ka)+i\lambda H_m^{(1)}(ka)}
\end{equation}
are the scattering coefficients of the disc.  A parameter scan gives, for
example,
$\min_m\Imag\kappa_m\approx-0.68$ at $k=1$, $a=29.4$, $\lambda=20$; since
$\kappa_m$ carries the factor $k^{-1/2}$, only the dimensionless combination
$\sqrt k\,\Imag\kappa_m$ is determined by $(ka,\lambda)$ alone.
By contrast,
$\Imag(e^{-i\pi/4}\kappa_m)=\sqrt{8\pi/k}\,(-\Real a_m)\ge0$ for every
passive disc, in agreement with \eqref{eq:optical}.  In three dimensions the
standard outgoing normalization has no corresponding phase rotation.
\end{remark}

Define the exact positive compact operator
\begin{equation}\label{eq:Fhash-exact}
 F_\#=|\Real\widetilde F|+\Imag\widetilde F.
\end{equation}
By Lemma~\ref{lem:optical}, $|\Imag\widetilde F|=\Imag\widetilde F$, so
\begin{equation}\label{eq:exact-equal-stabilized}
 F_\#=|\Real\widetilde F|+|\Imag\widetilde F|.
\end{equation}
For noisy empirical data the imaginary part need not remain positive, and we
therefore use the stabilized empirical definition
\begin{equation}\label{eq:empirical-sym}
 \widehat F_\#
 =|\Real(e^{-i\pi/4}\widehat F)|
 +|\Imag(e^{-i\pi/4}\widehat F)|.
\end{equation}
The second absolute value in \eqref{eq:empirical-sym} is a noise stabilization,
not an additional modification imposed by the exact impedance problem.

\subsection{Range identity}
The geometric role of $G$ is the following.

\begin{lemma}[Point-source range characterization]\label{lem:G-range}
For the sampling function \eqref{eq:test-function},
\begin{equation}\label{eq:G-geometry}
 z\in D\quad\Longleftrightarrow\quad\varphi_z\in\ran(G).
\end{equation}
\end{lemma}

\begin{proof}
Let $\Phi(x,z)=\tfrac i4H_0^{(1)}(k|x-z|)$ be the outgoing fundamental
solution.  In the physical normalization \eqref{eq:far-asymptotic} its far
field is $\gamma_2\varphi_z$ with $\gamma_2=e^{i\pi/4}/\sqrt{8\pi k}\ne0$.
If $z\in D$, then $\Phi(\cdot,z)$ is a radiating Helmholtz solution in a
neighborhood of $\overline{\R^2\setminus D}$ and smooth up to $\Gamma$, so
its impedance datum
$f_z=\partial_\nu\Phi(\cdot,z)+ik\lambda\,\Phi(\cdot,z)|_\Gamma$ lies in
$H^{-1/2}(\Gamma)$ and $Gf_z=\gamma_2\varphi_z$; hence
$\varphi_z=G(\gamma_2^{-1}f_z)\in\ran(G)$.

Conversely, suppose $\varphi_z=Gf$ for some $z\notin D$, and let $w_f$ be
the radiating solution with impedance datum $f$, so that $w_f$ has far
field $\varphi_z$ (cf.\ \cite[Chapter~2]{KG08}).  The field
$\gamma_2^{-1}\Phi(\cdot,z)$ has the same far field, so Rellich's lemma and
real-analytic unique continuation give
$w_f=\gamma_2^{-1}\Phi(\cdot,z)$ on the unbounded component $\Omega_0$ of
$\R^2\setminus(\overline D\cup\{z\})$.  If $z\in\R^2\setminus\overline D$,
choose a disc $B_\rho(z)\subset\R^2\setminus\overline D$.  Then
$w_f\in H^1(B_\rho(z))$, whereas
$|\nabla_x\Phi(x,z)|\asymp|x-z|^{-1}$ is not square integrable near $z$ in
two dimensions --- a contradiction.  If $z\in\Gamma$, apply the same
comparison on $B_\rho(z)\cap\Omega_0$, where the exterior solution $w_f$ is
$H^1$ up to the boundary, while the gradient of $\Phi(\cdot,z)$
again fails to be square integrable on that half-neighborhood, since a
fixed cone at $z$ contained in $\Omega_0$ does.  Hence $z\in D$.
\end{proof}

\begin{theorem}[Impedance factorization range identity]\label{thm:range}
Let $D$ and $\lambda$ satisfy the assumptions above.  Then
\begin{equation}\label{eq:range-identity}
 \ran(F_\#^{1/2})=\ran(G),
\end{equation}
and therefore
\begin{equation}\label{eq:point-range}
 z\in D\quad\Longleftrightarrow\quad
 \varphi_z\in\ran(F_\#^{1/2}).
\end{equation}
\end{theorem}

\begin{proof}
For the local impedance law, the boundary space in \cite{CCH13} is
$V(\Gamma)=L^2(\Gamma)$, and the compact embedding required by its generalized
factorization theorem is
$H^{1/2}(\Gamma)\Subset L^2(\Gamma)$.  Lemma~\ref{lem:noeig} verifies that
$k^2$ is not an interior eigenvalue associated with $Z_k$ or with $Z_k^*$.
Therefore
\cite[Theorem~4.11 of the arXiv version]{CCH13} gives
\[
 \ran(F_{{\rm st},\#}^{1/2})=\ran(G_{\rm st}),
 \qquad
 F_{{\rm st},\#}=|\Real F_{\rm st}|+\Imag F_{\rm st}.
\]
By \eqref{eq:normalization-transfer} and \eqref{eq:rotated},
\[
 \widetilde F=\frac{1}{\sqrt{8\pi k}}F_{\rm st},
 \qquad
 F_\#=\frac{1}{\sqrt{8\pi k}}F_{{\rm st},\#},
 \qquad
 \ran(G)=\ran(G_{\rm st}).
\]
Multiplication by a positive scalar does not change the range of a square
root, so \eqref{eq:range-identity} follows.  Lemma~\ref{lem:G-range} then gives
\eqref{eq:point-range}.
\end{proof}

\begin{remark}[Nonnormality]\label{rem:nonnormal}
For nonsymmetric impedance obstacles, $F$ need not be normal, so
Theorem~\ref{thm:range} cannot be replaced by a normal-operator formula
without additional hypotheses; nor does non-self-adjointness of the middle
operator force nonnormality of $F$ --- the absorbing disc is Fourier
diagonal and hence normal.
\end{remark}

\subsection{Picard criterion and regularized field}
The operator $F_\#$ is compact, positive, and injective.  Positivity follows
from \eqref{eq:Fhash-exact}, and injectivity follows from
$F_\#\succeq\Imag\widetilde F$ together with \eqref{eq:strict-ImF}.  Let
$(\gamma_j,\psi_j)_{j\ge1}$ be an orthonormal eigensystem with
$\gamma_j>0$ and $\gamma_j\downarrow0$.  The Picard criterion gives
\begin{equation}\label{eq:picard}
 z\in D
 \quad\Longleftrightarrow\quad
 W(z):=\sum_{j=1}^\infty
 \frac{|\ip{\varphi_z}{\psi_j}|^2}{\gamma_j}<\infty.
\end{equation}
For $\alpha>0$, define
\begin{align}
 W_\alpha(z)
 &:=\ip{\varphi_z}{F_\#(F_\#+\alpha I)^{-2}\varphi_z}
 \label{eq:Walpha}\\
 &=\sum_{j=1}^\infty
 \frac{\gamma_j}{(\gamma_j+\alpha)^2}
 |\ip{\varphi_z}{\psi_j}|^2.
\end{align}
The scalar filter increases pointwise to $1/t$ as $\alpha\downarrow0$, so
\begin{equation}\label{eq:Wlimit}
 W_\alpha(z)\uparrow W(z),
 \qquad
 W(z)<\infty\Longleftrightarrow z\in D,
 \qquad
 W(z)=+\infty\Longleftrightarrow z\notin D.
\end{equation}
Boundary points are included in the divergent case because $D$ is open.

The raw reciprocal $1/W_\alpha$ is useful for imaging but is not uniformly
bounded a priori across an arbitrary random family.  We therefore use a
bounded reciprocal as the primary uncertainty-quantification field.

\begin{definition}[Bounded factorization indicator]\label{def:bounded-indicator}
For $B>0$, let
\begin{equation}\label{eq:hB}
 h_B(w)=
 \begin{cases}
 B,&0\le w\le B^{-1},\\
 w^{-1},&w>B^{-1},
 \end{cases}
 \qquad
 I_{\alpha,B}(z)=h_B(W_\alpha(z)).
\end{equation}
Then $0\le I_{\alpha,B}\le B$ and $h_B$ is globally $B^2$-Lipschitz on
$[0,\infty)$.
\end{definition}

The limiting bounded indicator is
\begin{equation}\label{eq:IB-limit}
 I_B(z)=
 \begin{cases}
 h_B(W(z)),&z\in D,\\
 0,&z\notin D.
 \end{cases}
\end{equation}
Pointwise, $I_{\alpha,B}(z)\to I_B(z)$ as $\alpha\downarrow0$.

\section{Pathwise random factorization and measurability}\label{sec:pathwise}

For $y\in Y$, denote by $F_y$, $G_y$, $F_{\#,y}$,
$W_\alpha(\cdot;y)$, and $I_{\alpha,B}(\cdot;y)$ the objects associated with
$D_y$.  The continuity argument concerns the physical exterior problem and
does not depend on the boundary-integral representation later used for
computation.

Choose $R>0$ so that every deformation is the identity near
$\partial B_R$ and every $\overline D_y$ lies strictly inside $B_R$.  Set
\[
 \widehat\Omega=B_R\setminus\overline{\widehat D},
 \qquad
 J_y=D\Psi_y,\qquad
 j_y=|\det J_y|,
 \qquad
 \mathbb A_y=j_yJ_y^{-1}J_y^{-T}.
\]
If $\widehat\tau$ denotes a unit tangent on $\partial\widehat D$, write
\[
 J_{\Gamma,y}=|J_y\widehat\tau|
\]
for the surface Jacobian.  Let
$\mathcal T_{k,R}:H^{1/2}(\partial B_R)\to H^{-1/2}(\partial B_R)$ be the
exact outgoing Dirichlet-to-Neumann map.

\begin{proposition}[Uniform forward continuity]\label{prop:forward-continuity}
Under Assumption~\ref{ass:parametric},
\begin{equation}\label{eq:F-continuity}
 Y\ni y\longmapsto F_y\in\mathcal K(\LL)
\end{equation}
is continuous in operator norm.  Moreover, the pulled-back exterior solution
operators satisfy
\begin{equation}\label{eq:uniform-forward-inverse}
 \sup_{y\in Y}\norm{\mathcal A_y^{-1}}_
 {\mathcal L(H^1(\widehat\Omega)^*,H^1(\widehat\Omega))}<\infty,
\end{equation}
where $\mathcal A_y$ is the variational operator defined below.  In
particular, $\sup_{y\in Y}\opnorm{F_y}<\infty$.
\end{proposition}

\begin{proof}
Let $w_{y,g}$ be the scattered field generated by the Herglotz incident field
$v_g$, and set $\widehat w_{y,g}=w_{y,g}\circ\Psi_y$ on
$\widehat\Omega$.  Since $\Psi_y$ is the identity near $\partial B_R$, the
exact radiation condition is imposed on the same outer boundary for every
$y$.  The pulled-back variational form is
\begin{align}
 a_y(U,V)
 &=
 \int_{\widehat\Omega}
 \left(\mathbb A_y\nabla U\cdot\nabla\overline V
 -k^2j_yU\overline V\right)\,\dd x
 \label{eq:pulled-form}\\
 &\quad
 -ik\int_{\partial\widehat D}
 \widehat\lambda\,J_{\Gamma,y}U\overline V\,\dd s
 -\bigl\langle
 \mathcal T_{k,R}\gamma_RU,\gamma_RV
 \bigr\rangle_{H^{-1/2},H^{1/2}}.
 \nonumber
\end{align}
The right-hand side is the bounded functional
\begin{equation}\label{eq:pulled-rhs}
 \ell_{y,g}(V)
 =
 \int_{\partial\widehat D}
 \left[
 \bigl(\partial_{\nu_y}v_g\bigr)\circ\Psi_y
 +ik\widehat\lambda\,(v_g\circ\Psi_y)
 \right]\overline V\,J_{\Gamma,y}\,\dd s.
\end{equation}
Here the normal on $\partial\widehat D$ is the pullback of the obstacle normal
pointing into the exterior; the annular outward normal is its negative.  With
this convention the solution satisfies
$a_y(\widehat w_{y,g},V)=\ell_{y,g}(V)$ for all
$V\in H^1(\widehat\Omega)$.

Assumption~\ref{ass:parametric} implies uniform ellipticity of
$\mathbb A_y$ and
\[
 y\longmapsto
 \mathbb A_y,\ j_y,\ J_{\Gamma,y}
\]
continuously in the relevant $L^\infty$ norms.  Notice that the material
impedance value
$\lambda_y\circ\Psi_y=\widehat\lambda$ is fixed, whereas the pulled-back
Robin form still contains the geometry-dependent factor
$J_{\Gamma,y}$.  The Herglotz trace maps and their normal derivatives depend
continuously on $y$, so
\[
 y\longmapsto\mathcal A_y
 \in\mathcal L\bigl(H^1(\widehat\Omega),H^1(\widehat\Omega)^*\bigr),
 \qquad
 y\longmapsto\ell_{y,\cdot}
 \in\mathcal L\bigl(\LL,H^1(\widehat\Omega)^*\bigr)
\]
are norm-continuous.

The G{\aa}rding inequality for the uniformly elliptic volume part, the compact
embedding $H^1(\widehat\Omega)\Subset L^2(\widehat\Omega)$, and the standard
properties of the exact outgoing Dirichlet-to-Neumann map show that every
$\mathcal A_y$ is Fredholm of index zero.  If
$\mathcal A_yU=0$, pushing $U$ forward gives an outgoing solution of the
homogeneous exterior impedance problem; uniqueness gives $U=0$.  Hence every
$\mathcal A_y$ is an isomorphism.  Norm-continuity of inversion on the open
group of bounded isomorphisms yields continuity of
$y\mapsto\mathcal A_y^{-1}$, and compactness of $Y$ gives
\eqref{eq:uniform-forward-inverse}.

Finally, the far field is a fixed bounded linear function of the Cauchy data
on $\partial B_R$, because the deformation is the identity there.  Composing
this fixed trace-to-pattern map with
$\mathcal A_y^{-1}\ell_{y,\cdot}$ proves operator-norm continuity of
$F_y$.  Its compactness follows either from the smoothness of the far-field
kernel or from compactness of the Herglotz trace followed by the bounded
solution and far-field maps.
\end{proof}

\begin{theorem}[Pathwise range identity and joint measurability]\label{thm:pathwise}
Let $Y(\omega)$ be the random parameter in
Assumption~\ref{ass:parametric}.  Then, for every $\omega$,
\begin{equation}\label{eq:pathwise-range}
 z\in D_{Y(\omega)}
 \quad\Longleftrightarrow\quad
 \varphi_z\in\ran(F_{\#,Y(\omega)}^{1/2}).
\end{equation}
For every $\alpha>0$, the maps
\[
 (z,\omega)\longmapsto W_\alpha(z;\omega),
 \qquad
 (z,\omega)\longmapsto I_{\alpha,B}(z;\omega)
\]
are jointly measurable and continuous in $z$ for fixed $\omega$.  The
extended-real field
\begin{equation}\label{eq:exact-random-W}
 W(z;\omega)=\lim_{n\to\infty}W_{1/n}(z;\omega)
\end{equation}
is jointly measurable and is finite exactly on $D_{Y(\omega)}$.
\end{theorem}

\begin{proof}
The pathwise range identity is Theorem~\ref{thm:range}.  By
Proposition~\ref{prop:forward-continuity}, $\omega\mapsto F_{Y(\omega)}$ is
strongly measurable in the compact-operator norm.  The maps
$F\mapsto\Real F$, $F\mapsto\Imag F$, and $A\mapsto|A|$ are norm-continuous
on bounded self-adjoint sets, although the last map is not globally
operator-Lipschitz.  Hence $\omega\mapsto F_{\#,Y(\omega)}$ is measurable.
For fixed $\alpha$, continuous functional calculus makes
$A\mapsto A(A+\alpha I)^{-2}$ norm-continuous on the positive cone.  Since
$z\mapsto\varphi_z$ is continuous in $\LL$, the regularized field is a
Carath\'eodory function --- measurable in $\omega$ for each $z$ and
continuous in $z$ for each $\omega$ --- and is therefore jointly
measurable.  The limit
\eqref{eq:exact-random-W} is monotone and therefore jointly measurable as an
extended-real function.  Its finiteness set follows from
\eqref{eq:pathwise-range} and the Picard criterion.
\end{proof}

\begin{remark}[Pathwise imaging versus operator averaging]\label{rem:pathwise-mean}
The theorem applies to each $F_{Y(\omega)}$.  In general, $\E[F_{Y}]$ is
not the far-field operator of any physical obstacle and need not admit the
same factorization or range characterization.  The statistically meaningful
operation in the present model is therefore
$F_Y\mapsto I_{\alpha,B}(\cdot;Y)\mapsto
(\E I_{\alpha,B},\Var I_{\alpha,B})$, not factorization of $\E F_Y$.
\end{remark}

\section{Quadrature-weighted finite-sample theory}\label{sec:discrete}

\subsection{Discrete Hilbert-space identification}
Let $\{E_{j,M}\}_{j=1}^M$ be a partition of $\Sone$ with positive arc-length
weights $w_{j,M}=|E_{j,M}|$, and set
\[
 e_{j,M}=w_{j,M}^{-1/2}\mathbf 1_{E_{j,M}}.
\]
The map
\begin{equation}\label{eq:JM}
 J_M:\C^M\longrightarrow\LL,
 \qquad
 J_Mc=\sum_{j=1}^Mc_je_{j,M},
\end{equation}
is an isometry onto the piecewise-constant subspace.  Let
$\widehat x_{j,M}\in E_{j,M}$ and define the sampled far-field kernel matrix
\[
 K_M=\bigl[u^\infty(\widehat x_{i,M},\widehat x_{j,M})\bigr]_{i,j=1}^M,
 \qquad
 W_M=\operatorname{diag}(w_{1,M},\ldots,w_{M,M}).
\]
The quadrature-weighted Nystr\"om approximation under the Euclidean
identification and the compatible sampling vector are
\begin{equation}\label{eq:weighted-matrix}
 A_M=W_M^{1/2}K_MW_M^{1/2},
 \qquad
 b_{z,M}=W_M^{1/2}
 \bigl(e^{-ik\widehat x_{j,M}\cdot z}\bigr)_{j=1}^M.
\end{equation}
Write $\widetilde A_M=J_MA_MJ_M^*$, extended by zero on
$\ran(J_M)^\perp$, and $\widetilde b_{z,M}=J_Mb_{z,M}$.

\begin{assumption}[Quadrature and test approximation]\label{ass:quadrature}
For a compact sampling window $K\subset\R^2$, there are deterministic
sequences $\eta_M,\zeta_M\downarrow0$ such that
\begin{equation}\label{eq:eta-zeta}
 \opnorm{\widetilde A_M-F}\le\eta_M,
 \qquad
 \sup_{z\in K}\norm{\widetilde b_{z,M}-\varphi_z}_{\LL}\le\zeta_M.
\end{equation}
\end{assumption}
For a continuous far-field kernel and shape-regular angular partitions,
\eqref{eq:eta-zeta} follows from compactness of $F$, convergence of the
piecewise-constant projections, and quadrature consistency.

For any square matrix or bounded operator $A$, define
\begin{equation}\label{eq:S-def}
 \cS(A)=|\Real(e^{-i\pi/4}A)|+|\Imag(e^{-i\pi/4}A)|,
 \qquad
 \cQ_\alpha(A)=\cS(A)(\cS(A)+\alpha I)^{-2},
\end{equation}
the rotation implementing the two-dimensional passivity convention of
Lemma~\ref{lem:optical}.  For exact passive data,
$\cS(F)=\Fh$ by \eqref{eq:exact-equal-stabilized}.  Since the rotation is a
unimodular scalar, it changes no operator norm, and all perturbation
estimates below are insensitive to it.  The discrete ideal field is
\begin{equation}\label{eq:WM}
 W_{\alpha,M}(z)=b_{z,M}^*\cQ_\alpha(A_M)b_{z,M}.
\end{equation}

\begin{remark}[Uniform full-circle sampling]\label{rem:raw-scaling}
If $w_{j,M}=2\pi/M$, then $A_M=(2\pi/M)K_M$ and
$b_{z,M}=\sqrt{2\pi/M}\,(e^{-ik\widehat x_j\cdot z})_{j=1}^M$.  A
computation using the raw matrix $K_M$ is algebraically equivalent only when
\begin{equation}\label{eq:alpha-scaling}
 \alpha_{\rm weighted}=\frac{2\pi}{M}\alpha_{\rm raw};
\end{equation}
with this rescaling the value of $W_{\alpha,M}$ is unchanged, while without
it the regularization parameter and operator norm are direction-count
dependent.
\end{remark}

\subsection{Statistical matrix model}
For a fixed geometry realization, suppose that the multistatic matrix is
acquired $N$ times independently.  We call one such acquisition a
\emph{measurement frame} --- the term refers to a repeated recording of the
full matrix, and carries no connection with frames in the sense of
overcomplete systems.  The $N$ frames produce matrices
$K_M^{(1)},\ldots,K_M^{(N)}$, each with mean
$\E K_M^{(\ell)}=K_M$.  Set
\begin{equation}\label{eq:empirical-A}
 \widehat K_{M,N}=\frac1N\sum_{\ell=1}^NK_M^{(\ell)},
 \qquad
 \widehat A_{M,N}=W_M^{1/2}\widehat K_{M,N}W_M^{1/2}.
\end{equation}
Let
$X_{\ell,M}=W_M^{1/2}(K_M^{(\ell)}-K_M)W_M^{1/2}$.  A concrete sufficient
condition for a deterministic confidence radius is the standard
matrix-Bernstein hypothesis
\begin{equation}\label{eq:bernstein-hyp}
\begin{aligned}
 &\E X_{\ell,M}=0,
 \qquad
 \opnorm{X_{\ell,M}}\le L_M\ \text{a.s.},\\
 &v_M=\max\left\{
 \opnorm{\E(X_{\ell,M}X_{\ell,M}^*)},\;
 \opnorm{\E(X_{\ell,M}^*X_{\ell,M})}
 \right\},
\end{aligned}
\end{equation}
the last quantity being the matrix variance parameter.
Then a rectangular matrix Bernstein inequality \cite{Tropp15} gives, for a universal $C>0$,
\begin{equation}\label{eq:bernstein-radius}
 \Prob\left\{
 \opnorm{\widehat A_{M,N}-A_M}>b_{M,N}(\delta)
 \right\}\le\delta,
\end{equation}
where one may take
\begin{equation}\label{eq:bMN}
 b_{M,N}(\delta)=C\left[
 \sqrt{\frac{v_M\log(2M/\delta)}{N}}
 +\frac{L_M\log(2M/\delta)}{N}
 \right].
\end{equation}
Sub-Gaussian or sub-exponential frame models can be substituted, with the
almost-sure bound in \eqref{eq:bernstein-hyp} replaced by a tail (Orlicz-norm)
hypothesis and a confidence radius of the same order; the
fixed-relative-Frobenius experiment of Section~\ref{sec:numerics} is a
sensitivity check, not a literal realization of \eqref{eq:bernstein-hyp}.
The dimension dependence is determined by the concentration theorem and the
noise covariance (compare \cite{KL17}) and cannot in general be replaced by
an aperture-independent effective rank.

\subsection{Two perturbation lemmas}
The operator absolute value is norm-continuous \cite{Kato73,AP11} but is not globally
operator-Lipschitz.  The sharp dimension-free modulus on a bounded set has a
double-logarithmic correction.

\begin{lemma}[Modulus of continuity of the symmetrization]\label{lem:modulus}
There is a universal constant $C>0$ such that, for every $\Lambda>0$, whenever
$\max\{\opnorm A,\opnorm B\}\le\Lambda$,
\begin{equation}\label{eq:rho-bound}
 \opnorm{\cS(A)-\cS(B)}
 \le \rho_\Lambda(\opnorm{A-B}),
\end{equation}
where
\begin{equation}\label{eq:rho}
 \rho_\Lambda(t)=C\,\tau
 \log\!\left(2+\log\!\frac{2\Lambda}{\tau}\right),
 \qquad \tau:=\min\{t,2\Lambda\},
 \qquad \rho_\Lambda(0)=0.
\end{equation}
The clipping at $2\Lambda$ costs nothing, because
$\opnorm{A-B}\le2\Lambda$ always holds under the stated norm bounds; it
makes $\rho_\Lambda$ well defined and nondecreasing on all of $[0,\infty)$,
so that it may be evaluated at any deterministic upper bound for the error,
even one exceeding $2\Lambda$.
\end{lemma}

\begin{proof}
Write $\widetilde A=e^{-i\pi/4}A$ and $\widetilde B=e^{-i\pi/4}B$; the
rotation changes neither $\opnorm{A-B}$ nor the norm bound $\Lambda$.  Apply
the sharp operator-modulus estimate of Aleksandrov and Peller
\cite{AP11} to the self-adjoint pairs
$(\Real\widetilde A,\Real\widetilde B)$ and
$(\Imag\widetilde A,\Imag\widetilde B)$.  Because
$\opnorm{\Real\widetilde A-\Real\widetilde B}$ and
$\opnorm{\Imag\widetilde A-\Imag\widetilde B}$ do not exceed
$\opnorm{A-B}$, while each of the
four self-adjoint operators has norm at most $\Lambda$, summing the two
absolute-value estimates proves \eqref{eq:rho-bound} after changing the
universal constant.  For the monotonicity, differentiation gives
\[
 \rho_\Lambda'(t)
 =C\log\Bigl(2+\log\frac{2\Lambda}{t}\Bigr)
 -\frac{C}{2+\log(2\Lambda/t)}
 \ \ge\ C\Bigl(\log2-\tfrac12\Bigr)>0
 \qquad(0<t\le2\Lambda),
\]
since $\log(2\Lambda/t)\ge0$ there.
\end{proof}

\begin{lemma}[Filtered resolvent stability]\label{lem:filter}
If $C,D\succeq0$ and $\alpha>0$, then
\begin{equation}\label{eq:filter-stability}
 \opnorm{C(C+\alpha I)^{-2}-D(D+\alpha I)^{-2}}
 \le\frac{3}{\alpha^2}\opnorm{C-D}.
\end{equation}
Moreover,
\begin{equation}\label{eq:filter-norm}
 \opnorm{C(C+\alpha I)^{-2}}\le\frac{1}{4\alpha}.
\end{equation}
\end{lemma}

\begin{proof}
Since
$t(t+\alpha)^{-2}=(t+\alpha)^{-1}-\alpha(t+\alpha)^{-2}$,
the first resolvent identity bounds the first difference by
$\alpha^{-2}\opnorm{C-D}$.  Expanding the difference of the squared resolvents
gives two further terms, each bounded after multiplication by $\alpha$; their
total is at most $2\alpha^{-2}\opnorm{C-D}$.  The scalar maximum of
$t/(t+\alpha)^2$ is $1/(4\alpha)$.
\end{proof}

The double logarithm in Lemma~\ref{lem:modulus} is a global worst-case statement.  In
a fixed finite-dimensional active subspace with a spectral gap from zero, the
usual local Lipschitz behavior is recovered.

\begin{proposition}[Local gap-dependent absolute-value bound]\label{prop:gap}
Let $H$ and $\widetilde H$ be self-adjoint matrices with
$\max\{\opnorm H,\opnorm{\widetilde H}\}\le\Lambda$ and
\[
 \dist(0,\sigma(H)),\ \dist(0,\sigma(\widetilde H))\ge g>0.
\]
Then
\begin{equation}\label{eq:gap-bound}
 \opnorm{|H|-|\widetilde H|}
 \le C\left(1+\frac{\Lambda}{g}\right)
 \opnorm{H-\widetilde H}
\end{equation}
for a universal constant $C$.  Applying this separately to real and imaginary
parts gives a local Lipschitz estimate for $\cS$.
\end{proposition}

\begin{proof}
Use $|H|=H\sgn(H)$.  Apply the corresponding identity to
$\widetilde H$.  The Davis--Kahan/Riesz-projection estimate gives
\[
 \opnorm{\sgn(H)-\sgn(\widetilde H)}
 \le Cg^{-1}\opnorm{H-\widetilde H}.
\]
The result follows by adding and subtracting
$\widetilde H\sgn(H)$.
\end{proof}

\subsection{Fully discrete deterministic error bound}
Let
\begin{equation}\label{eq:empirical-W}
 \widehat W_{\alpha,M,N}(z)
 =b_{z,M}^*\cQ_\alpha(\widehat A_{M,N})b_{z,M}.
\end{equation}
Set
\[
 R_K=\sup_{z\in K}\norm{\varphi_z}_{\LL}=\sqrt{2\pi},
 \qquad
 e_{M,N}=\eta_M+s_{M,N},
\]
where $s_{M,N}$ is any deterministic bound for
$\opnorm{\widehat A_{M,N}-A_M}$.

The norm radius used in the absolute-value modulus must be deterministic.
For a fixed obstacle one may take any $\Lambda$ larger than
$\opnorm F+\sup_M\opnorm{\widetilde A_M}$ and then enlarge it by the prescribed
data-error radius.  For the random family, Proposition~\ref{prop:forward-continuity}
and the uniform discretization assumption introduced in
Section~\ref{sec:end-to-end} provide a common choice independent of the
realization.

\begin{theorem}[Fully discrete perturbation bound]\label{thm:discrete-bound}
Assume that the exact, lifted discrete, and empirical operators have norm at
most $\Lambda$.  On the event
$\opnorm{\widehat A_{M,N}-A_M}\le s_{M,N}$,
\begin{align}
 \sup_{z\in K}
 |\widehat W_{\alpha,M,N}(z)-W_\alpha(z)|
 &\le
 \frac{3(R_K+\zeta_M)^2}{\alpha^2}
 \rho_\Lambda(e_{M,N})
 \label{eq:discrete-bound}\\
 &\quad+
 \frac{(2R_K+\zeta_M)\zeta_M}{4\alpha}.
\end{align}
\end{theorem}

\begin{proof}
Lift the empirical matrix and test vector by $J_M$.  The operator difference
from $F$ is at most $e_{M,N}$, and $\rho_\Lambda$ is nondecreasing, so
Lemmas~\ref{lem:modulus} and
\ref{lem:filter} bound the filtered-operator difference by
$3\alpha^{-2}\rho_\Lambda(e_{M,N})$.  Multiplication by the squared test-vector
norm gives the first term.  For the test-vector replacement, use
\[
 |\ip{Qx}{x}-\ip{Qy}{y}|
 \le\opnorm Q\norm{x-y}(\norm x+\norm y)
\]
and \eqref{eq:filter-norm}.
\end{proof}

Since $h_B$ is $B^2$-Lipschitz, Theorem~\ref{thm:discrete-bound} immediately
gives
\begin{equation}\label{eq:I-error}
 \sup_{z\in K}
 |\widehat I_{\alpha,B,M,N}(z)-I_{\alpha,B}(z)|
 \le B^2 E_{\alpha,M,N},
\end{equation}
where $E_{\alpha,M,N}$ denotes the right side of
\eqref{eq:discrete-bound}.

\subsection{Probability modes and joint consistency}
A convergence statement for a sequence of experiments must use deterministic
confidence radii and a specified probability mode.

\begin{theorem}[Joint $M$--$N$--$\alpha$ consistency]\label{thm:joint-consistency}
Let $(M_n,N_n)$ be a sequence.  Suppose
\[
 \sup_n\max\{\opnorm F,\opnorm{\widetilde A_{M_n}}\}\le\Lambda_0<\infty
\]
and deterministic numbers $b_n\downarrow0$ and $\delta_n\in(0,1)$ satisfy
\begin{equation}\label{eq:confidence-sequence}
 \Prob\left\{
 \opnorm{\widehat A_{M_n,N_n}-A_{M_n}}>b_n
 \right\}\le\delta_n.
\end{equation}
For all sufficiently large $n$, set
$\Lambda=\Lambda_0+1$; on the confidence event below,
$\opnorm{\widehat A_{M_n,N_n}}\le\Lambda$.  Choose
$\alpha_n\downarrow0$ so that
\begin{equation}\label{eq:joint-coupling}
 \frac{\rho_\Lambda(\eta_{M_n}+b_n)}{\alpha_n^2}\longrightarrow0,
 \qquad
 \frac{\zeta_{M_n}}{\alpha_n}\longrightarrow0.
\end{equation}
Then, for every fixed $z$,
\begin{equation}\label{eq:joint-pointwise}
 \widehat W_{\alpha_n,M_n,N_n}(z)
 \longrightarrow
 \begin{cases}
 W(z)<\infty,&z\in D,\\
 +\infty,&z\notin D,
 \end{cases}
\end{equation}
in probability if $\delta_n\to0$.  If
$\sum_n\delta_n<\infty$, the convergence is almost sure.  The bounded empirical
indicator converges correspondingly to $I_B(z)$.
\end{theorem}

\begin{proof}
On the event in \eqref{eq:confidence-sequence},
Theorem~\ref{thm:discrete-bound} and \eqref{eq:joint-coupling} give
$\widehat W_{\alpha_n,M_n,N_n}(z)-W_{\alpha_n}(z)\to0$.  Equation
\eqref{eq:Wlimit} gives the deterministic limit.  Failure probabilities tending
to zero yield convergence in probability.  If they are summable,
Borel--Cantelli shows that all but finitely many error events fail to occur,
which gives almost-sure convergence.  Apply continuity and boundedness of
$h_B$ for the indicator statement.
\end{proof}

\begin{remark}[Admissible versus balanced schedules]\label{rem:schedule}
A schedule such as $\alpha\asymp\rho_\Lambda(\beta)^{1/3}$ may satisfy a
perturbation coupling, but no optimality statement is available without a
bias model, and the classical H\"older source conditions that would supply
one are empty for far-field imaging.  For the impedance disc they fail at
every noncentral interior point, while logarithmic source conditions hold
for every exponent and yield only logarithmic bias rates.  These questions,
together with the associated sharp disc asymptotics, are studied
separately in \cite{PartII}; here we record only that any admissible
schedule suffices for the consistency statements below.
\end{remark}

\section{Uncertainty quantification and confidence bands}\label{sec:uq}

\subsection{Monte Carlo mean and variance on a finite grid}
Let $Y_1,\ldots,Y_Q$ be independent copies of the geometry parameter.  For a
fixed $\alpha$ and $B$, set
\begin{equation}\label{eq:Xq}
 X_q(z)=I_{\alpha,B}(z;Y_q)\in[0,B].
\end{equation}
The empirical mean and unbiased sample variance are
\begin{equation}\label{eq:empirical-moments}
 \overline X_Q(z)=\frac1Q\sum_{q=1}^QX_q(z),
 \qquad
 S_Q^2(z)=\frac{1}{Q-1}\sum_{q=1}^Q
 \bigl(X_q(z)-\overline X_Q(z)\bigr)^2.
\end{equation}
Equivalently,
\begin{equation}\label{eq:pairwise-variance}
 S_Q^2(z)=\frac{1}{Q(Q-1)}
 \sum_{1\le i<j\le Q}(X_i(z)-X_j(z))^2.
\end{equation}

\begin{theorem}[Mean and variance consistency]\label{thm:mc}
For every fixed $z$,
\[
 \overline X_Q(z)\longrightarrow\E X_1(z),
 \qquad
 S_Q^2(z)\longrightarrow\Var X_1(z)
\]
almost surely as $Q\to\infty$.  Moreover,
\begin{equation}\label{eq:mean-mse}
 \E\left|\overline X_Q(z)-\E X_1(z)\right|^2
 =\frac{\Var X_1(z)}{Q}\le\frac{B^2}{4Q}.
\end{equation}
\end{theorem}

\begin{proof}
The variables are bounded, so the strong law applies to $X_q$ and $X_q^2$.
The variance identity follows by substitution.  Equation
\eqref{eq:mean-mse} is the variance of an average together with Popoviciu's
bound $\Var X\le B^2/4$.
\end{proof}

\begin{theorem}[Simultaneous finite-grid concentration]\label{thm:grid-concentration}
Let $\mathcal Z$ be a finite sampling grid with $J=|\mathcal Z|$.  With
probability at least $1-\delta$,
\begin{equation}\label{eq:mean-concentration}
 \max_{z\in\mathcal Z}
 |\overline X_Q(z)-\E X_1(z)|
 \le B\sqrt{\frac{\log(2J/\delta)}{2Q}}.
\end{equation}
With probability at least $1-\delta$,
\begin{equation}\label{eq:variance-concentration}
 \max_{z\in\mathcal Z}
 |S_Q^2(z)-\Var X_1(z)|
 \le B^2\sqrt{\frac{\log(2J/\delta)}{2Q}}.
\end{equation}
\end{theorem}

\begin{proof}
Apply Hoeffding's inequality at each grid point and take a union bound to obtain
\eqref{eq:mean-concentration}.  For the variance, use
\eqref{eq:pairwise-variance}.  Replacing one observation changes $S_Q^2(z)$
by at most $B^2/Q$.  McDiarmid's bounded-differences inequality therefore gives
$2\exp(-2Qt^2/B^4)$ at a fixed point.  A union bound yields
\eqref{eq:variance-concentration}.
\end{proof}

\subsection{Continuum confidence bands}
A finite-grid bound becomes a uniform bound on a compact sampling
window once the spatial modulus of the indicator is made explicit.  Let
$N(K,r)$ be the least cardinality of an $r$-net of a compact set
$K\subset\R^2$.

\begin{lemma}[Spatial Lipschitz bound]\label{lem:spatial-lipschitz}
For every admissible realization and all $z,z'\in\R^2$,
\begin{align}
 |W_\alpha(z)-W_\alpha(z')|
 &\le \frac{\pi k}{\alpha}|z-z'|,
 \label{eq:W-spatial-lipschitz}\\
 |I_{\alpha,B}(z)-I_{\alpha,B}(z')|
 &\le \frac{\pi kB^2}{\alpha}|z-z'|.
 \label{eq:I-spatial-lipschitz}
\end{align}
The same bounds hold for the quadrature-weighted matrix field when
$\norm{b_{z,M}}^2=\sum_jw_{j,M}=2\pi$ and
$|b_{z,M,j}-b_{z',M,j}|\le k\sqrt{w_{j,M}}|z-z'|$.
\end{lemma}

\begin{proof}
The elementary inequality $|e^{it}-e^{is}|\le|t-s|$ gives
\[
 \norm{\varphi_z-\varphi_{z'}}_{\LL}
 \le k\sqrt{2\pi}\,|z-z'|,
 \qquad
 \norm{\varphi_z}_{\LL}=\sqrt{2\pi}.
\]
Lemma~\ref{lem:filter} gives
$\norm{F_\#(F_\#+\alpha I)^{-2}}\le(4\alpha)^{-1}$.  Hence
\[
 |W_\alpha(z)-W_\alpha(z')|
 \le\frac{1}{4\alpha}
 \norm{\varphi_z-\varphi_{z'}}
 \bigl(\norm{\varphi_z}+\norm{\varphi_{z'}}\bigr)
 \le\frac{\pi k}{\alpha}|z-z'|.
\]
Composition with the $B^2$-Lipschitz map $h_B$ proves
\eqref{eq:I-spatial-lipschitz}.  The matrix statement is identical.
\end{proof}

\begin{theorem}[Uniform mean-field confidence band]\label{thm:continuum-band}
Let $K\subset\R^2$ be compact, let $r>0$, and put
$m(z)=\E I_{\alpha,B}(z;Y)$.  With probability at least $1-\delta$,
\begin{equation}\label{eq:continuum-band}
 \sup_{z\in K}|\overline X_Q(z)-m(z)|
 \le
 B\sqrt{\frac{\log\bigl(2N(K,r)/\delta\bigr)}{2Q}}
 +\frac{2\pi kB^2}{\alpha}r.
\end{equation}
If empirical pathwise indicators satisfy
$\max_q\sup_{z\in K}|\widehat X_q(z)-X_q(z)|\le d$, then the right side of
\eqref{eq:continuum-band} is increased by $d$.
\end{theorem}

\begin{proof}
Apply Theorem~\ref{thm:grid-concentration} to an $r$-net of $K$.  Both the
sample mean and its expectation are $\pi kB^2/\alpha$-Lipschitz by
Lemma~\ref{lem:spatial-lipschitz}.  Moving from an arbitrary point to its
nearest net point therefore costs at most twice this Lipschitz constant times
$r$.  The empirical-pathwise statement follows from the triangle inequality.
\end{proof}

\begin{remark}[Target and conservativeness of the continuum band]
\label{rem:net-radius}
The target in Theorem~\ref{thm:continuum-band} is
$m_\alpha(z)=\E I_{\alpha,B}(z;Y)$, the mean \emph{regularized indicator}.
It is not the occupancy probability $p(z)=\Prob\{z\in D_Y\}$ and is not,
without an additional level-set argument, a confidence region for the random
obstacle.  For a planar window of diameter $D_K$ one may use
$N(K,r)\le C_K(1+D_K/r)^2$; choosing $r=\alpha/(2\pi kB\sqrt Q)$ turns the
interpolation term into $B/\sqrt Q$ while the covering term grows only
logarithmically in $1/r$, so the band tends to zero whenever
$\alpha=\alpha_Q$ with $\log(1/\alpha_Q)=o(Q)$.  The global Lipschitz
constant scales as $\alpha^{-1}$, so the bound is conservative at very
small regularization.  Localized moduli, empirical-process chaining, or the
occupancy construction in Section~\ref{sec:occupancy} can yield more directly
interpretable sets.
\end{remark}

\subsection{Combined measurement and geometry sampling error}
Let $Y_1,\ldots,Y_Q$ be independent geometry draws.  Conditional on
$Y_q$, suppose the $N$ measurement frames are independent, have the same
conditional mean matrix, and are independent across different geometry
groups.  Let
$\widehat X_q(z)=\widehat I_{\alpha,B,M,N}(z;Y_q)$.  Allocate an
operator-estimation failure probability $\delta_{\rm op}/Q$ to each path.
The number $d$ below is a deterministic \emph{uniform-in-geometry} indicator
error bound obtained from \eqref{eq:I-error}; thus the quadrature, norm, and
matrix-concentration constants entering it must be bounded uniformly over
$y\in Y$.

\begin{theorem}[Combined confidence bound]\label{thm:combined}
Assume
\begin{equation}\label{eq:pathwise-d}
 \Prob\left\{
 \max_{z\in\mathcal Z}|\widehat X_q(z)-X_q(z)|>d
 \right\}\le\frac{\delta_{\rm op}}{Q}
\end{equation}
for every $q$.  Let $\widehat{\overline X}_Q$ and $\widehat S_Q^2$ be the mean
and unbiased variance formed from the empirical pathwise indicators.  With
probability at least $1-\delta_{\rm op}-\delta_{\rm mc}$,
\begin{align}
 \max_{z\in\mathcal Z}
 |\widehat{\overline X}_Q(z)-\E X_1(z)|
 &\le d+B\sqrt{\frac{\log(4J/\delta_{\rm mc})}{2Q}},
 \label{eq:combined-mean}\\
 \max_{z\in\mathcal Z}
 |\widehat S_Q^2(z)-\Var X_1(z)|
 &\le2Bd+B^2\sqrt{\frac{\log(4J/\delta_{\rm mc})}{2Q}}.
 \label{eq:combined-var}
\end{align}
\end{theorem}

\begin{proof}
A union bound over the $Q$ paths gives simultaneous pathwise accuracy with
failure probability at most $\delta_{\rm op}$.  On that event, the two sample
means differ by at most $d$.  For the variances, every pairwise squared
difference changes by at most $4Bd$ because both exact and empirical indicators
lie in $[0,B]$.  Formula \eqref{eq:pairwise-variance} then gives
$|\widehat S_Q^2-S_Q^2|\le2Bd$.  Apply
Theorem~\ref{thm:grid-concentration} to the exact pathwise sample with failure
probability $\delta_{\rm mc}$, splitting it between mean and variance.
\end{proof}

\begin{remark}[Three sample sizes]\label{rem:budgets}
The direction count $M$, frame count $N$, and geometry-realization count $Q$
play different roles.  Increasing $M$ reduces quadrature and test-vector error;
increasing $N$ reduces the pathwise data-operator error; increasing $Q$ reduces
Monte Carlo fluctuation of the geometry push-forward.  Theorem~\ref{thm:combined}
keeps these three sample sizes separate.  Repeated independent reconstructions at the
same geometry introduce a fourth design parameter, used next to separate
measurement and geometry variability.
\end{remark}

\subsection{Joint stochastic consistency}
\label{sec:end-to-end}

The fixed-geometry estimates become a stochastic-geometry theorem once their
constants are uniform over the admissible family.

\begin{assumption}[Uniform finite-data approximation over the geometry family]
\label{ass:uniform-discrete}
For a compact sampling window $K$, there are deterministic
$\eta_M,\zeta_M\downarrow0$ such that
\begin{equation}\label{eq:uniform-discrete}
 \sup_{y\in Y}\opnorm{\widetilde A_{M,y}-F_y}\le\eta_M,
 \qquad
 \sup_{y\in Y}\sup_{z\in K}
 \norm{\widetilde b_{z,M}-\varphi_z}_{\LL}\le\zeta_M.
\end{equation}
There is a deterministic $\Lambda<\infty$ bounding the exact and ideal
discrete operators uniformly in $y$ and $M$.  Moreover, for every
$\delta\in(0,1)$ there is a deterministic radius $b_{M,N}(\delta)$ such that
\begin{equation}\label{eq:uniform-conditional-concentration}
 \sup_{y\in Y}
 \Prob\!\left(
 \opnorm{\widehat A_{M,N}(y)-A_M(y)}
 >b_{M,N}(\delta)
 \,\middle|\,Y=y
 \right)
 \le\delta.
\end{equation}
\end{assumption}

The first norm bound follows from Proposition~\ref{prop:forward-continuity}
for a uniformly convergent Nystr\"om family; the second is geometry-independent
for a fixed angular rule.  Condition
\eqref{eq:uniform-conditional-concentration} is the precise uniformity needed
when the same confidence radius is used for every random path.

For a finite grid $\mathcal Z\subset K$, let
\[
 \widehat m_n(z)=\frac1{Q_n}\sum_{q=1}^{Q_n}
 \widehat I_{\alpha_n,B,M_n,N_n}(z;Y_{n,q}),
\]
and let $\widehat v_n(z)$ be the corresponding unbiased sample variance.  The
arrays $\{Y_{n,q}\}_{q=1}^{Q_n}$ are independent copies of $Y$; conditional
measurement frames are independent within and across groups.

\begin{theorem}[Joint $M$--$N$--$Q$--$\alpha$ consistency]
\label{thm:end-to-end}
Under Assumption~\ref{ass:uniform-discrete}, let
$M_n,N_n,Q_n\to\infty$, $\alpha_n\downarrow0$, and choose
$\delta_{{\rm op},n}\downarrow0$.  Put
\[
 b_n=b_{M_n,N_n}\!\left(\frac{\delta_{{\rm op},n}}{Q_n}\right)
\]
and suppose that, with the enlarged radius $\Lambda_*:=\Lambda+1$,
\begin{equation}\label{eq:end-to-end-coupling}
 \frac{\rho_{\Lambda_*}(\eta_{M_n}+b_n)}{\alpha_n^2}\to0,
 \qquad
 \frac{\zeta_{M_n}}{\alpha_n}\to0.
\end{equation}
Then, uniformly on the finite grid $\mathcal Z$,
\begin{align}
 \widehat m_n(z)&\longrightarrow
 m_B(z):=\E I_B(z;Y), \label{eq:end-to-end-mean}\\
 \widehat v_n(z)&\longrightarrow
 v_B(z):=\Var(I_B(z;Y)) \label{eq:end-to-end-var}
\end{align}
in probability.  If
\[
 \sum_{n=1}^\infty\delta_{{\rm op},n}<\infty
 \qquad\text{and}\qquad
 \frac{Q_n}{\log n}\longrightarrow\infty,
\]
then both convergences hold almost surely.

\end{theorem}

\begin{proof}
Let
$X_{n,q}(z)=I_{\alpha_n,B}(z;Y_{n,q})$.
Assumption~\ref{ass:uniform-discrete} bounds only the exact and ideal
operators by $\Lambda$; since \eqref{eq:end-to-end-coupling} forces
$b_n\to0$, we have $b_n\le1$ for all large $n$, and on each confidence event
$\opnorm{\widehat A_{M_n,N_n}(y)}\le\Lambda+1=\Lambda_*$, so
Theorem~\ref{thm:discrete-bound} applies with norm radius $\Lambda_*$.  By
that theorem, \eqref{eq:I-error}, and a union bound over $q$, with
probability at least $1-\delta_{{\rm op},n}$,
\[
 \max_{q\le Q_n}\max_{z\in\mathcal Z}
 |\widehat I_{\alpha_n,B,M_n,N_n}(z;Y_{n,q})
 -X_{n,q}(z)|
 \le d_n,
\]
where $d_n\to0$ by \eqref{eq:end-to-end-coupling}.  On this event, empirical
means differ from the ideal means by at most $d_n$, and empirical sample
variances differ from the ideal sample variances by at most $2Bd_n$.

For the ideal triangular array, Hoeffding's inequality and a union bound give,
for every $t>0$,
\[
 \Prob\!\left(
 \max_{z\in\mathcal Z}
 \left|
 \frac1{Q_n}\sum_{q=1}^{Q_n}X_{n,q}(z)
 -\E X_{n,1}(z)
 \right|>t
 \right)
 \le2|\mathcal Z|e^{-2Q_nt^2/B^2}.
\]
The corresponding McDiarmid bound for the sample variance is
$2|\mathcal Z|e^{-2Q_nt^2/B^4}$.  Finally,
$I_{\alpha_n,B}(z;Y)\to I_B(z;Y)$ almost surely and is bounded by $B$.
Dominated convergence therefore gives, uniformly on the finite set
$\mathcal Z$,
\[
 \E X_{n,1}(z)\to m_B(z),
 \qquad
 \Var X_{n,1}(z)\to v_B(z).
\]
Combining the three errors proves convergence in probability.  Under the
summability hypotheses, Borel--Cantelli applied first to the operator events
and then to the concentration bounds for a countable sequence of tolerances
gives almost-sure convergence.
\end{proof}

\subsection{Separating geometry and measurement uncertainty}
The variance of a reconstructed field mixes physical variability of the
obstacle with variability induced by repeated measurements.  A nested design
separates the two without linearizing the imaging map.

Let the groups
$\bigl(Y_q,X_{q1},\ldots,X_{qR}\bigr)$, $q=1,\ldots,Q$, be independent and
identically distributed.  Conditional on $Y_q$, let
$X_{q\ell}(z)\in[0,B]$, $\ell=1,\ldots,R$, be conditionally independent and
identically distributed reconstructed indicators; one replicate may itself
be formed from an $N$-frame averaged operator.  Define
\[
 \mu_z(y)=\E[X_{q\ell}(z)\mid Y_q=y],
 \qquad
 \sigma_z^2(y)=\Var(X_{q\ell}(z)\mid Y_q=y),
\]
and
\begin{equation}\label{eq:nested-components}
 V_{\rm geom}(z)=\Var_Y\mu_z(Y),
 \qquad
 V_{\rm meas}(z)=\E_Y\sigma_z^2(Y).
\end{equation}
The law of total variance gives
\begin{equation}\label{eq:total-variance}
 \Var(X_{q\ell}(z))=V_{\rm geom}(z)+V_{\rm meas}(z).
\end{equation}

For $Q,R\ge2$, set
\begin{align*}
 \overline X_q(z)&=\frac1R\sum_{\ell=1}^R X_{q\ell}(z),\\
 S_{q,{\rm w}}^2(z)&=\frac{1}{R-1}\sum_{\ell=1}^R
 \bigl(X_{q\ell}(z)-\overline X_q(z)\bigr)^2,\\
 \overline{\overline X}(z)&=\frac1Q\sum_{q=1}^Q\overline X_q(z),\\
 S_{\rm b}^2(z)&=\frac{1}{Q-1}\sum_{q=1}^Q
 \bigl(\overline X_q(z)-\overline{\overline X}(z)\bigr)^2.
\end{align*}

\begin{theorem}[Unbiased nested variance decomposition]\label{thm:nested-variance}
The estimators
\begin{equation}\label{eq:nested-estimators}
 \widehat V_{\rm meas}(z)=\frac1Q\sum_{q=1}^Q S_{q,{\rm w}}^2(z),
 \qquad
 \widehat V_{\rm geom}(z)
 =S_{\rm b}^2(z)-\frac{1}{R}\widehat V_{\rm meas}(z)
\end{equation}
are unbiased for $V_{\rm meas}(z)$ and $V_{\rm geom}(z)$, respectively.
For fixed $R$, both are strongly consistent as $Q\to\infty$.  Consequently
\[
 \widehat V_{\rm total}
 :=\widehat V_{\rm geom}+\widehat V_{\rm meas}
 =S_{\rm b}^2+\left(1-\frac1R\right)\widehat V_{\rm meas}
\]
is unbiased and strongly consistent for the total variance in
\eqref{eq:total-variance}.
\end{theorem}

\begin{proof}
Conditionally on $Y_q$, the within-geometry sample variance is unbiased, so
$\E S_{q,{\rm w}}^2=V_{\rm meas}$.  Moreover,
\[
 \Var(\overline X_q)
 =\Var_Y\mu_z(Y)+\frac1R\E_Y\sigma_z^2(Y)
 =V_{\rm geom}+\frac1R V_{\rm meas}.
\]
The sample variance $S_{\rm b}^2$ is unbiased for this quantity.  Subtracting
$R^{-1}\widehat V_{\rm meas}$ proves unbiasedness.  Strong consistency follows
from the strong law applied to the bounded variables
$\overline X_q$, $\overline X_q^2$, and $S_{q,{\rm w}}^2$.
\end{proof}

\begin{remark}[Nonnegative displays versus unbiased estimation]
The unbiased estimator $\widehat V_{\rm geom}$ can be slightly negative at
finite $Q$ when the true geometry component is small.  The untruncated value
should be reported for inferential calculations.  The clipped field
$\max\{0,\widehat V_{\rm geom}\}$ may be used for visualization, but it is no
longer unbiased.
\end{remark}

\begin{theorem}[Finite-grid concentration for the nested variance components]
\label{thm:nested-concentration}
Let $\mathcal Z$ be a finite grid with $J$ points and assume the independent
group design above.  With probability at least $1-\delta$, simultaneously for
all $z\in\mathcal Z$,
\begin{align}
 |\widehat V_{\rm meas}(z)-V_{\rm meas}(z)|
 &\le
 B^2\sqrt{\frac{\log(4J/\delta)}{2Q}},
 \label{eq:nested-meas-concentration}\\
 |\widehat V_{\rm geom}(z)-V_{\rm geom}(z)|
 &\le
 B^2\sqrt{\frac{2\log(4J/\delta)}{Q}}.
 \label{eq:nested-geom-concentration}
\end{align}
\end{theorem}

\begin{proof}
Regard each group $\mathcal G_q=(Y_q,X_{q1},\ldots,X_{qR})$ as one
independent observation.  Since $0\le S_{q,{\rm w}}^2\le B^2$, replacing one
group changes $\widehat V_{\rm meas}$ by at most $B^2/Q$, and McDiarmid's
inequality gives
$\Prob\{|\widehat V_{\rm meas}-V_{\rm meas}|>t\}\le2e^{-2Qt^2/B^4}$.  By the
pairwise representation
$S_{\rm b}^2=\tfrac{1}{Q(Q-1)}\sum_{i<j}(\overline X_i-\overline X_j)^2$,
replacing one group changes $S_{\rm b}^2$ by at most $B^2/Q$ and
$R^{-1}\widehat V_{\rm meas}$ by at most $B^2/(RQ)$, hence
$\widehat V_{\rm geom}$ by at most $2B^2/Q$, and McDiarmid gives
$\Prob\{|\widehat V_{\rm geom}-V_{\rm geom}|>t\}\le2e^{-Qt^2/(2B^4)}$.
Unbiasedness (Theorem~\ref{thm:nested-variance}) and a union bound over the
grid and the two components conclude.
\end{proof}

\begin{remark}[Cost-optimal number of repeats per geometry]
\label{rem:optimal-repeats}
If one new geometry costs $c_{\rm g}$ and each conditionally independent
reconstruction at a geometry costs $c_{\rm m}$, then under the total cost
$C=Q(c_{\rm g}+Rc_{\rm m})$ the grand mean has variance
$(c_{\rm g}+Rc_{\rm m})\bigl(V_{\rm geom}+V_{\rm meas}/R\bigr)/C$, minimized
over $R$ at
$R_*(z)=\sqrt{c_{\rm g}V_{\rm meas}(z)/(c_{\rm m}V_{\rm geom}(z))}$.  This
is the optimum of the continuous relaxation in $R$; only integer repeat
counts are feasible, and the admissible optimum is found by comparing the
objective at the two integers adjacent to $R_*(z)$.  Nor is it directly
usable, since both components are unknown and depend on
$z$; in practice pilot estimates are inserted, and a single global $R$ is
obtained by applying the rule to integrated or maximum-risk variance
components.  In the edge cases the relaxation degenerates.  Taking $V_{\rm geom}=0$
gives $R_*=\infty$, while $V_{\rm meas}=0$ places the stationary point at
$R=0$, which is infeasible, so the objective is then minimized at the
smallest allowed repeat count --- $R=1$ for estimating the mean, and $R=2$
whenever the variance components themselves are to be estimated.
\end{remark}

\subsection{Interpretation, occupancy probability, and variance leakage}
\label{sec:occupancy}
Let $I_0(z)$ denote a nominal-geometry indicator and
$m(z)=\E I_{\alpha,B}(z;Y)$.  The exact decomposition
\begin{equation}\label{eq:bias-variance}
 \E|I_{\alpha,B}(z;Y)-I_0(z)|^2
 =\Var(I_{\alpha,B}(z;Y))+|m(z)-I_0(z)|^2
\end{equation}
shows that the variance is a path-to-path sensitivity diagnostic, not total
reconstruction error.

For the membership field
\begin{equation}\label{eq:membership}
 M(z;Y)=\mathbf1_{\{z\in D_Y\}},
 \qquad p(z)=\Prob\{z\in D_Y\},
\end{equation}
one has the exact Bernoulli identity
\begin{equation}\label{eq:membership-var}
 \Var M(z;Y)=p(z)(1-p(z)).
\end{equation}
Thus membership variance vanishes outside the membership-uncertainty set
$\{0<p<1\}$.  This exact zero-set statement does not transfer automatically
to a continuous regularized indicator.

\begin{definition}[Thresholded occupancy estimator]
For $\alpha,\tau>0$, define
\begin{equation}\label{eq:threshold-membership}
 M_{\alpha,\tau}(z;y)
 =\mathbf1_{\{W_\alpha(z;y)\le\tau\}}.
\end{equation}
This field is used only for occupancy inference; the continuous bounded
indicator remains preferable for visualization and Lipschitz perturbation
bounds.
\end{definition}

\begin{theorem}[Diagonal recovery of occupancy probability]
\label{thm:occupancy}
Let $\mathcal Z$ be a finite sampling grid.  There exist deterministic
sequences $\tau_n\uparrow\infty$ and $\alpha_n\downarrow0$ such that
\begin{equation}\label{eq:occupancy-classification}
 \max_{z\in\mathcal Z}
 \Prob\!\left\{
 M_{\alpha_n,\tau_n}(z;Y)
 \ne\mathbf1_{\{z\in D_Y\}}
 \right\}\longrightarrow0.
\end{equation}
Consequently, uniformly on $\mathcal Z$,
\begin{align}
 \E M_{\alpha_n,\tau_n}(z;Y)&\longrightarrow
 p(z):=\Prob\{z\in D_Y\}, \label{eq:occupancy-mean}\\
 \Var M_{\alpha_n,\tau_n}(z;Y)&\longrightarrow
 p(z)(1-p(z)). \label{eq:occupancy-var}
\end{align}

Suppose an empirical field satisfies
$\Prob\{\max_{z\in\mathcal Z}
|\widehat W_n(z;Y)-W_{\alpha_n}(z;Y)|>d_n\}\le r_n$ with $d_n,r_n\to0$
(a pointwise bound at level $r_n'$ yields this with
$r_n=|\mathcal Z|\,r_n'$ by a union bound).  If
\begin{equation}\label{eq:occupancy-margin}
 \max_{z\in\mathcal Z}
 \Prob\{|W_{\alpha_n}(z;Y)-\tau_n|\le d_n\}\longrightarrow0,
\end{equation}
then the empirical classifier
$\widehat M_n=\mathbf1_{\{\widehat W_n\le\tau_n\}}$ satisfies the same
conclusions.  The sequences $(\tau_n,\alpha_n)$ are deterministic but may
depend on the law of $Y$ and on $\mathcal Z$; the statement is an existence
theorem, not a prescription for choosing a threshold from data.
\end{theorem}

\begin{proof}
On $\{z\in D_Y\}$ the exact Picard field $W(z;Y)$ is finite, so
$\Prob\{z\in D_Y,\ W(z;Y)>\tau\}\to0$ as $\tau\to\infty$; for fixed $\tau$,
monotone convergence $W_\alpha(z;Y)\uparrow W(z;Y)$ gives
$\Prob\{z\notin D_Y,\ W_\alpha(z;Y)\le\tau\}\to0$ as $\alpha\downarrow0$.
Choose $\tau_n\uparrow\infty$ so that the first probability is at most $1/n$
on $\mathcal Z$, then $\alpha_n\downarrow0$ so that the second is at most
$1/n$ on $\mathcal Z$; this proves \eqref{eq:occupancy-classification}.
Both fields being Bernoulli, $L^1$ convergence gives
\eqref{eq:occupancy-mean}, and convergence of the means gives
\eqref{eq:occupancy-var}.  For empirical fields, the ideal and empirical
classifiers can disagree only on the operator-error event or when
$|W_{\alpha_n}-\tau_n|\le d_n$; \eqref{eq:occupancy-margin} and $r_n\to0$
finish the proof.
\end{proof}

\begin{remark}[Occupancy confidence bands]
For $Q$ independent paths the empirical mean of $M_{\alpha_n,\tau_n}$ is a
Bernoulli proportion, so Hoeffding or exact binomial intervals give
finite-sample bands whose target, in the diagonal limit of
Theorem~\ref{thm:occupancy}, converges to the physical occupancy $p(z)$ ---
distinct from the mean-indicator band of Theorem~\ref{thm:continuum-band}.
\end{remark}

\begin{proposition}[Variance leakage bound]\label{prop:variance-leakage}
For every $z$, define the essential pathwise oscillation with respect to
the law $P_Y$ of $Y$,
\[
 \omega_{\alpha,B}(z)=\operatorname*{ess\,sup}_{P_Y}I_{\alpha,B}(z;\cdot)
 -\operatorname*{ess\,inf}_{P_Y}I_{\alpha,B}(z;\cdot).
\]
Then
\begin{equation}\label{eq:osc-var}
 \Var(I_{\alpha,B}(z;Y))\le\frac{\omega_{\alpha,B}(z)^2}{4}.
\end{equation}
In particular, if the indicator oscillation is at most $\epsilon$ outside an
enlarged boundary-uncertainty band, then the variance leakage there is at most
$\epsilon^2/4$.
\end{proposition}

\begin{proof}
Equation \eqref{eq:osc-var} is Popoviciu's inequality applied to the essential
range at the fixed point $z$.
\end{proof}

\section{Numerical experiments}\label{sec:numerics}

All experiments use full-circle observation and incidence directions.  The
phrase ``finite angular sampling'' is used deliberately, since missing angular
sectors would constitute a different limited-aperture inverse problem and
would require restriction/data-recovery operators or a separate factorization
analysis \cite{LS17}.  Unless a spectral norm is explicitly
indicated, relative matrix errors use the Frobenius norm.  Three test
configurations are used: a random three-lobed star, the kite of
\cite{CK19}, and a two-component scene pairing the kite with a tilted
ellipse.  The reference parametrizations and placements, all experiment
parameters, and the validation of the discretizations --- analytic disc far
fields and layer-operator eigenvalues, reciprocity, boundary-node and
two-component refinement, and the far-separation interaction defect --- are
collected in Appendices~\ref{app:validation} and~\ref{app:params}.

\subsection{Solver, resonance qualification, and weighted imaging}
Two boundary-integral representations are used.  The first represents the
scattered field by a single-layer potential $u^s=S_k\phi$.  The exterior
normal trace gives
\begin{equation}\label{eq:bie}
 \left(-\frac12I+K_k'+ik\lambda S_k\right)\phi
 =-\left(\partial_\nu u^i+ik\lambda u^i\right).
\end{equation}
For a smooth periodic parametrization, the logarithmic singularity of $S_k$ is
integrated by the trigonometric Kress rule \cite{Kress14}, while smooth
off-diagonal multiple-scattering blocks use the periodic trapezoidal rule.
The far field is evaluated by the same periodic quadrature.

Physical uniqueness of the exterior impedance problem does not by itself make
the indirect equation \eqref{eq:bie} resonance-free.  The following standard
qualification separates the physical problem from the chosen representation.

\begin{proposition}[Representation-induced resonance condition]
\label{prop:single-layer-resonance}
On a smooth boundary, the operator
\[
 -\frac12I+K_k'+ikM_\lambda S_k:
 H^{-1/2}(\Gamma)\longrightarrow H^{-1/2}(\Gamma)
\]
is Fredholm of index zero.  It is injective, and hence invertible, if and only
if $k^2$ is not an interior Dirichlet eigenvalue of $D$.  At every interior
Dirichlet eigenvalue it has a nontrivial nullspace, even though the exterior
impedance problem remains uniquely solvable.
\end{proposition}

\begin{proof}
On a smooth boundary, $K_k'$ and $S_k$ are compact relative to the identity on
the natural trace space, so the operator is Fredholm of index zero.  If a
homogeneous density produces a scattered field, the boundary equation implies
that this field satisfies the homogeneous exterior impedance condition and is
therefore zero outside $D$.  The single-layer potential is continuous across
$\Gamma$, so its interior trace is zero.  If $k^2$ is not an interior
Dirichlet eigenvalue, the interior potential is zero and the jump in its normal
derivative gives $\phi=0$.  Conversely, let $u$ be a nonzero interior Dirichlet eigenfunction.  Its
normal derivative does not vanish identically, since otherwise $u$ would have
zero Cauchy data and vanish by Cauchy uniqueness; it is therefore a nonzero
single-layer density whose potential has zero
exterior trace and hence vanishes outside, and this density lies in the
nullspace of \eqref{eq:bie}.
\end{proof}

Thus the factorization theory itself is valid for every real $k>0$, whereas
the particular solver \eqref{eq:bie} is used only after conditioning has been
checked at the selected frequencies.  A resonance-free alternative starts
from the combined potential $u^s=(K_k-i\eta S_k)\phi$ with $\eta>0$.  On a
smooth boundary define
\begin{equation}\label{eq:Acf}
 \mathcal A^{\rm cf}_{D,k,\eta}
 :=
 N_k-i\eta\left(-\frac12I+K_k'\right)
 +ikM_\lambda\left(\frac12I+K_k-i\eta S_k\right)
 :
 H^{1/2}(\Gamma)\longrightarrow H^{-1/2}(\Gamma).
\end{equation}
The combined-field equation is
\begin{equation}\label{eq:combined-field-bie}
 \mathcal A^{\rm cf}_{D,k,\eta}\phi
 =-\left(\partial_\nu u^i+ik\lambda u^i\right).
\end{equation}

\begin{proposition}[All-frequency combined-field isomorphism]
\label{prop:combined-isomorphism}
Let $\Gamma$ be of class $C^\infty$.  For every $k,\eta>0$ and every positive real
$\lambda\in L^\infty(\Gamma)$ bounded away from zero,
$\mathcal A^{\rm cf}_{D,k,\eta}$ is a bounded Fredholm operator of index zero
from $H^{1/2}(\Gamma)$ to $H^{-1/2}(\Gamma)$ and is an isomorphism.
\end{proposition}

\begin{proof}
On a $C^\infty$ boundary the hypersingular operator $N_k$ is
an elliptic pseudodifferential operator of order one and is Fredholm of
index zero from $H^{1/2}(\Gamma)$ to $H^{-1/2}(\Gamma)$; the mapping and
Fredholm properties of the layer operators used in this argument are those
recorded in \cite{CCH13,Kress14}.  The remaining terms in \eqref{eq:Acf}
are compact between the same spaces.  The operators $-\tfrac12I+K_k'$ and
$\tfrac12I+K_k-i\eta S_k$ are bounded on $H^{1/2}(\Gamma)$, multiplication
by $\lambda\in L^\infty(\Gamma)$ is bounded on $L^2(\Gamma)$, and the
embeddings
$H^{1/2}(\Gamma)\Subset L^2(\Gamma)\hookrightarrow H^{-1/2}(\Gamma)$ are
compact and bounded, respectively; composing these maps exhibits every
non-hypersingular term of \eqref{eq:Acf} as a compact operator
$H^{1/2}(\Gamma)\to H^{-1/2}(\Gamma)$.  Hence
$\mathcal A^{\rm cf}_{D,k,\eta}$ is Fredholm of index zero.

Let $\phi$ solve the homogeneous equation and put
$v^\pm=(K_k-i\eta S_k)\phi$ in the exterior and interior, respectively.
The exterior trace satisfies the homogeneous physical impedance condition, so
$v^+=0$ by exterior uniqueness.  With the jump conventions used in
\eqref{eq:bie},
\[
 \gamma_0^-v=-\phi,
 \qquad
 \gamma_1^-v=-i\eta\phi,
\]
and therefore
\[
 \partial_\nu v^- -i\eta v^-=0\qquad\text{on }\Gamma.
\]
Green's identity, with $\eta>0$, gives
$v^-=0$ exactly as in Lemma~\ref{lem:noeig}.  Thus $\phi=0$.
Injectivity and index zero imply surjectivity.
\end{proof}

\begin{remark}[Boundary regularity]\label{rem:cf-regularity}
Only the Fredholm step above uses more than the standing $C^2$ regularity,
through the pseudodifferential description of $N_k$, and that step survives
at lower regularity.  On a Lipschitz boundary, the single-layer operator and
the normal derivative of the double-layer potential of a strongly elliptic
second-order operator --- here $\Delta+k^2$ --- both satisfy G\aa rding
inequalities in their energy norms \cite{Costabel88}; the inequality for the
latter gives $N_k:H^{1/2}(\Gamma)\to H^{-1/2}(\Gamma)$ Fredholm of index
zero.  The compactness of the remaining terms of \eqref{eq:Acf} and the
injectivity argument use only the jump relations and Green's identity, and
are unchanged.  A G\aa rding inequality, rather than plain coercivity, is
what is available here and what the argument needs.  In two dimensions the
single-layer operator fails to be positive definite on all of
$H^{-1/2}(\Gamma)$ when the logarithmic capacity of $\Gamma$ is not small,
a scaling effect that the compact remainder absorbs.
\end{remark}

\looseness=-1
The first-order operator in \eqref{eq:Acf} may be right-regularized by an
elliptic isomorphism $\mathcal R_{D,k}:H^{-1/2}(\Gamma)\to H^{1/2}(\Gamma)$
--- a coercive modified-Helmholtz single layer is standard --- producing a
second-kind equation, as in multifrequency impedance-obstacle computations
\cite{BorgesRachh22}.  Proposition~\ref{prop:combined-isomorphism} gives invertibility at every
parameter.  If in addition the pulled-back combined-field family is
norm-continuous in $(y,k,\eta)$ --- a property that the argument of
Proposition~\ref{prop:forward-continuity} establishes for the volume
variational form, but which we do not prove here for the boundary-operator
family --- then compactness of
$Y\times[k_{\min},k_{\max}]\times[\eta_{\min},\eta_{\max}]$ yields uniform
inverse bounds.  \emph{Provided} the regularized product-integration
family is collectively compact and uniformly consistent over this compact
parameter set, collectively compact approximation theory
\cite[Chapter~12]{Kress14} then yields uniform discrete stability and
far-field convergence.  Verifying collective compactness for a specific
regularizer $\mathcal R_{D,k}$ and quadrature family remains
open.  The computations reported here discretize
\eqref{eq:combined-field-bie} directly.  The hypersingular operator is
assembled through Maue's identity with the Kress rule and trigonometric
differentiation.  The resulting blocks reproduce the analytic disc
eigenvalues of $K_k$ and $N_k$ to near machine precision.  The far field
agrees with the analytic disc series --- including at interior Dirichlet
eigenvalues, where the single-layer system is numerically singular --- and
with the single-layer solver away from resonance.  The ensemble
reconstructions and the repeated-noise study reported below are generated
with this combined-field discretization; the monitored single-layer solver
is retained for validation and for the pathwise conditioning comparison in
Table~\ref{tab:conditioning-ensemble}.

Every sampled kernel $K_M$ and every sampling vector are converted to the
weighted forms \eqref{eq:weighted-matrix}.  We diagonalize
\[
 \cS(A_M)=|\Real(e^{-i\pi/4}A_M)|+|\Imag(e^{-i\pi/4}A_M)|,
\]
with the phase rotation applied before either Hermitian part is formed, and
compute the bounded indicator with $B=0.35$.  The ceiling is chosen so that
$1/B$ is comparable to typical interior values of $W_\alpha$; it clips deep
interior reciprocal values, suppressing spurious interior variance at the
cost of a flat saturated interior mean.  Under uniform sampling, the values
agree with a raw-matrix calculation only after the parameter rescaling
\eqref{eq:alpha-scaling}; the measured relative difference in a disc test is
$5.5\times10^{-15}$.

\subsection{Analytic disc reference and pathwise conditioning}
For a disc of radius $a$ the scattering coefficients $a_m$ are those of
\eqref{eq:disc-am}, and the exact far-field kernel is obtained by Fourier
summation.  In the
physical normalization \eqref{eq:far-asymptotic} the far-field operator is
Fourier diagonal, with eigenvalues and symmetrized eigenvalues
\begin{equation}\label{eq:disc-kappa}
 \kappa_m=2\pi\sqrt{\frac{2}{\pi k}}\,e^{-i\pi/4}a_m,
 \qquad
 \widetilde\kappa_m=e^{-i\pi/4}\kappa_m,
 \qquad
 \gamma_m=|\Real\widetilde\kappa_m|+\Imag\widetilde\kappa_m,
\end{equation}
the $\gamma_m$ being the eigenvalues of $\Fh$; passivity gives
$\Imag\widetilde\kappa_m=\sqrt{8\pi/k}\,(-\Real a_m)\ge0$, so the
second absolute value in $\Fh$ is inactive for the disc.  The boundary
operators have eigenvalues
\begin{equation}\label{eq:disc-layer-eigs}
 s_m=\frac{i\pi a}{2}J_m(ka)H_m^{(1)}(ka),
 \qquad
 \kappa_m'=-\frac12+
 \frac{i\pi ka}{2}J_m'(ka)H_m^{(1)}(ka).
\end{equation}
The $-1/2$ term in $\kappa_m'$ is essential; omitting it tests a different
jump convention.

Because of Proposition~\ref{prop:single-layer-resonance}, conditioning must
be examined separately from physical well-posedness.

\looseness=-1
Nominal systems alone do not characterize a random-domain computation, since
the interior Dirichlet spectrum and hence the conditioning of
\eqref{eq:bie} move with the geometry.  Recomputing every ensemble
path from its recorded seed gives Table~\ref{tab:conditioning-ensemble}.
These are Euclidean two-norm condition numbers of the actual dense Nystr\"om
matrices; they are not condition numbers of the continuum operator.  The two
formulations discretize operators of different orders between different
trace spaces, so their singular values carry different scalings and the
cross-formulation comparison below is indicative rather than invariant; what
it does establish is the presence or absence of near-singular systems within
each formulation.

\begin{table}[tbhp]
\tbl{Median, p90 (90th percentile), and maximum of the condition
number over the $Q$ paths of each random ensemble; $\sigma_{\min}$ is the
smallest singular value over all paths.  Single-layer
rows refer to the solver \eqref{eq:bie}; combined-field rows
re-assemble the resonance-free system \eqref{eq:combined-field-bie} with
$\eta=k$ on the identical paths.}
{\small
\setlength{\tabcolsep}{4pt}
\begin{tabular}{@{}lrrrrrr@{}}
\toprule
geometry / solver & $Q$ & median & p90 & maximum
& min.\ $\sigma_{\min}$ & worst seed\\
\midrule
star, single layer & 48 & 5.66 & 6.70 & 7.94
& $3.10\times10^{-1}$ & 3042\\
\quad combined field & 48 & 5.17 & 5.37 & 5.51
& $3.98$ & 3042\\
kite, single layer & 48 & 35.8 & 149 & 841
& $3.76\times10^{-3}$ & 11040\\
\quad combined field & 48 & 6.49 & 6.64 & 6.77
& $6.53$ & 11007\\
two components, single layer & 36 & 54.0 & 309 & $1.99\times10^3$
& $1.86\times10^{-3}$ & 22008\\
\quad combined field & 36 & 5.33 & 5.48 & 5.57
& $6.52$ & 22014\\
\bottomrule
\end{tabular}}
\label{tab:conditioning-ensemble}
\end{table}

\looseness=-1
The kite and two-component ensembles include paths whose single-layer
systems are substantially worse conditioned than the nominal ones, and the
observed worst-case conditioning can grow sharply as the ensemble is
enlarged.  Over an extended run of fresh random kite draws, the worst single-layer condition number rises
from $2.4\times10^{2}$ over the first $48$ draws to $1.1\times10^{5}$
over $2000$, while the combined-field condition number on the identical
draws never exceeds $6.8$ (Table~\ref{tab:conditioning-tail}).  We interpret this as random geometry sampling closer to the moving interior
Dirichlet spectrum as the ensemble grows; on that reading the validity of a
single-layer ensemble is a property of its particular seeds rather than of
the method.  We do not claim an asymptotic divergence, which would require
knowing that the geometry law charges every neighborhood of a resonant
parameter.  All reported ensembles are therefore
generated with the combined field \eqref{eq:combined-field-bie}, with
$\sigma_{\min}$ and the condition number recorded for every realization; the
single-layer rows of Table~\ref{tab:conditioning-ensemble} document the
representation that these seeds would otherwise have used.

\begin{table}[tbhp]
\tbl{Worst-case conditioning of a growing random-kite ensemble at
$k=7$ (fresh seeds, disjoint from those of the reported ensembles).  Running maximum
of $\operatorname{cond}_2$ for the two formulations on the identical
draws.}
{\begin{tabular}{lrrrrr}
\toprule
draws & 48 & 132 & 500 & 1000 & 2000\\
\midrule
single layer & 236 & $2.60\times10^{3}$ & $9.23\times10^{3}$
& $9.23\times10^{3}$ & $1.15\times10^{5}$\\
combined field & 6.70 & 6.72 & 6.76 & 6.77 & 6.79\\
\bottomrule
\end{tabular}}
\label{tab:conditioning-tail}
\end{table}

\subsection{Admissible geometry priors}\label{sec:priors-num}
The three-lobed star uses
\begin{equation}\label{eq:star-prior}
 r(t;y)=1.45+0.34\cos(3t)
 +0.10\sum_{n=1}^4\frac{\xi_n\cos(nt)+\eta_n\sin(nt)}{n},
 \qquad \xi_n,\eta_n\sim\mathcal U[-1,1].
\end{equation}
The modal amplitude bound gives
\[
 r(t;y)\ge1.45-0.34-0.10\sqrt2\sum_{n=1}^4\frac1n
 >0.815.
\]
Hence every realization is a regular simple star-shaped curve whose speed is
bounded below by its radius.  (The weaker
per-realization bound $a_0-\sum_n(|c_n|+|s_n|)\ge0.69$ was checked draw by
draw.)  The family is embedded in Assumption~\ref{ass:parametric} by an
explicit composition of radial near-identity maps.

\begin{lemma}[Explicit admissible extension of the star prior]\label{lem:star-map}
Write $r_0(t)=1.45+0.34\cos(3t)$ and
$\delta_y(t)=0.10\sum_{n=1}^4(\xi_n\cos nt+\eta_n\sin nt)/n$, so that the
realized radius in \eqref{eq:star-prior} is $r_y=r_0+\delta_y$.  Fix a
cut-off $\chi\in C^\infty([0,\infty))$ with $0\le\chi\le1$, $\chi=0$ on
$[0,0.55]\cup[4,\infty)$, $\chi=1$ on $[0.8,2.2]$, $|\chi'|\le9$ on
$[0.55,0.8]$, and $|\chi'|\le1.2$ on $[2.2,4]$.  In polar coordinates
$x=\rho\,e_\rho(t)$, $e_\rho(t)=(\cos t,\sin t)$, define
\[
 \Phi_y=\mathrm{Id}+V_y,
 \qquad
 V_y(x)=\tfrac16\,\chi(\rho)\,\delta_y(t)\,e_\rho(t),
 \qquad
 \Psi_y=\Phi_y^{\circ6}\quad(\text{six-fold composition}).
\]
Then $\sup_y\norm{DV_y}_{L^\infty}\le0.49$; each $\Phi_y$, and hence
$\Psi_y$, is a global bi-Lipschitz $C^2$ diffeomorphism equal to the identity
outside $B_4$; $\Psi_y$ maps the reference curve $\rho=r_0(t)$ onto the
realized curve $\rho=r_y(t)$ and the reference star domain onto the realized
one; and $y\mapsto\Psi_y$ is continuous in $C^2$.  In particular the star
prior satisfies Assumption~\ref{ass:parametric} with constants independent of
the realization.
\end{lemma}

\begin{proof}
The $n$th mode of $\delta_y$ has pointwise amplitude at most $0.10\sqrt2/n$
and differentiation in $t$ multiplies it by $n$, so
$|\delta_y|\le0.2947$, $|\delta_y'|\le0.5657$, and $|\delta_y''|\le1.4143$.
In the orthonormal frame $(e_\rho,e_t)$, with $e_t=(-\sin t,\cos t)$, the
derivative of $\chi(\rho)\delta_y(t)\,e_\rho$ is the upper-triangular matrix
$\bigl(\begin{smallmatrix}\chi'\delta_y&\chi\delta_y'/\rho\\
0&\chi\delta_y/\rho\end{smallmatrix}\bigr)$; bounding its spectral norm by
its Frobenius norm on the three regions $\rho\in[0.55,0.8]$, $[0.8,2.2]$,
and $[2.2,4]$ gives $2.90$, $0.80$, and $0.46$, respectively, so
$\sup_y\norm{DV_y}_{L^\infty}\le2.90/6\le0.49$.  The field vanishes for
$\rho\le0.55$ and outside $B_4$, and its second derivatives are bounded
uniformly in $y$ by the displayed mode bounds and the fixed cut-off.
Lemma~\ref{lem:diffeomorphism} with $q=0.49$ then makes each $\Phi_y$, hence
the composition $\Psi_y$, a global bi-Lipschitz $C^2$ diffeomorphism with
constants $(1-q)^6$ and $(1+q)^6$ uniform in $y$.  Along the ray of angle
$t$, $\Phi_y$ acts as the strictly increasing scalar map
$\rho\mapsto\rho+\tfrac16\chi(\rho)\delta_y(t)$.  Since
$r_0\in[1.11,1.79]$ and $r_y=r_0+\delta_y\in[0.8153,2.0847]\subset[0.8,2.2]$,
every intermediate radius $r_0+\tfrac j6\delta_y$ lies where $\chi\equiv1$;
the six iterates therefore move the boundary radius by six equal increments
$\delta_y/6$, so $\Psi_y(r_0(t)\,e_\rho(t))=r_y(t)\,e_\rho(t)$, and
monotonicity along rays maps the reference domain onto the realized one.
Finally, $y\mapsto V_y$ is affine, hence continuous, into $C^2$, and
composition preserves this; conditions (A1)--(A5) of
Assumption~\ref{ass:parametric} follow, the impedance being constant.
\end{proof}

A normal-displacement prior on the kite is not used, because such a prior can
leave the tubular neighborhood of a concave reference curve and self-intersect.
The kite and two-component configurations use the affine family
\begin{equation}\label{eq:affine-prior}
 A_y=I+
 \begin{pmatrix}
 0.05y_1&0.025y_2\\
 0.025y_3&0.05y_4
 \end{pmatrix},
 \qquad
 b_y=0.06(y_5,y_6)^T,
 \qquad y_j\sim\mathcal U[-1,1].
\end{equation}
Since
$\norm{A_y-I}_2\le\norm{A_y-I}_F\le0.0791$, one has
\begin{equation}\label{eq:affine-margin}
 \sigma_{\min}(A_y)\ge0.9209.
\end{equation}
Thus affine images cannot self-intersect; and since one draw moves the whole
scene --- kite and ellipse by the same map --- the component separation is at
least $0.9209$ times the nominal separation.  The smallest singular value
among $10^4$ validation draws was $0.9332$.  Remark~\ref{rem:affine-extension}
embeds this boundary family in the compactly supported domain-map construction.

\subsection{Pathwise and ensemble reconstructions}
\looseness=-1
For the star we use $k=4$, $M=64$, $96$ boundary nodes, $Q=48$, and
$\alpha=\frac{2\pi}{64}(2.0\times10^{-3})=1.9635\times10^{-4}$; for the kite
and two-component cases, $k=7$, $M=96$, $Q=48$ and $36$, respectively, and
$\alpha=\frac{2\pi}{96}(1.2\times10^{-3})=7.8540\times10^{-5}$.  For the two-component scene the kite component uses $80$ nodes and the
ellipse $56$, which the refinement study of Appendix~\ref{app:validation} resolves to
$2\times10^{-8}$.  All ensembles are generated with the combined-field
solver; the condition number of every solved system is recorded with the
realization, and all $132$ values lie below $6.8$
(Table~\ref{tab:conditioning-ensemble}).  All ensemble variances use the
unbiased denominator $Q-1$.

\looseness=-1
The Monte Carlo resolution of the displayed fields is measured by the
pointwise standard error $\sqrt{S_Q^2(z)/Q}$, whose maximum over the grid is
$5.1\times10^{-3}$ for the star ($2.3\%$ of the peak mean; $Q=48$, $J=8100$),
$5.7\times10^{-3}$ for the kite ($3.2\%$; $Q=48$, $J=7740$), and
$7.5\times10^{-3}$ for the two-component scene ($4.3\%$; $Q=36$, $J=9240$)
--- pointwise diagnostics, not simultaneous bands.  The distribution-free
finite-grid mean-band radius $B\sqrt{\log(2J/\delta)/(2Q)}$ of
Theorem~\ref{thm:grid-concentration} at $\delta=0.05$ is $0.127$--$0.148$
(the variance-band radii carry $B^2$ in place of $B$), because it protects roughly $10^4$ grid points at once;
the continuum band of Theorem~\ref{thm:continuum-band} adds a spatial
interpolation term and is more conservative still.  The gap between
pointwise and simultaneous uncertainty is reported rather than hidden.

\begin{figure}[tbhp]
\centering
\includegraphics[width=0.6\textwidth]{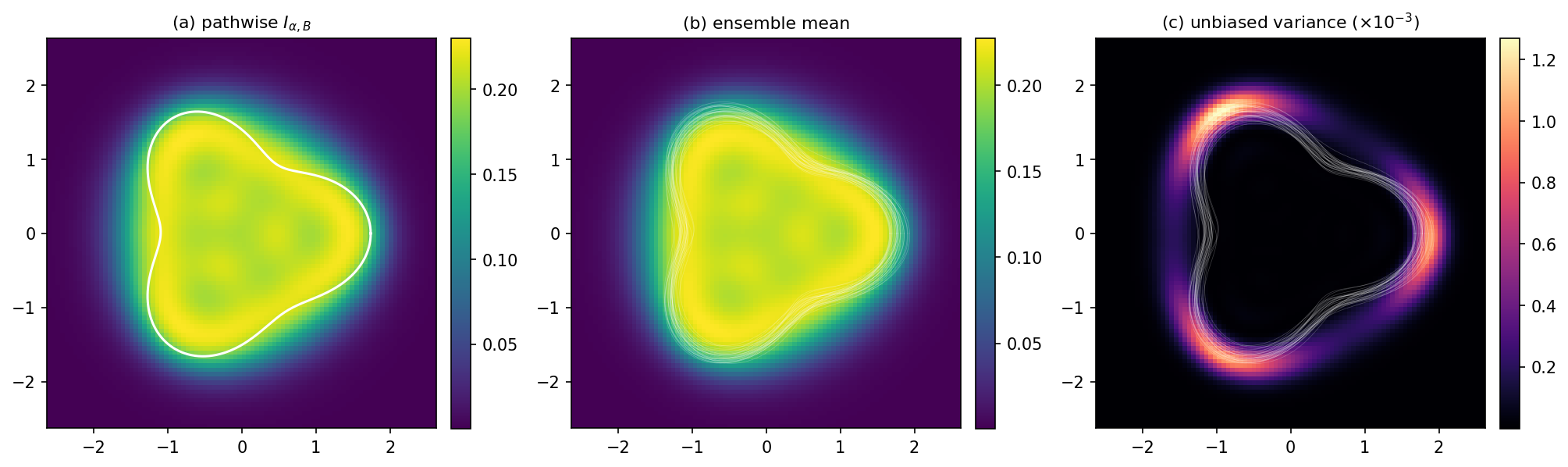}
\caption{Admissible random-star reconstruction.  White curves show the
pathwise boundary in panel (a) and a subset of ensemble boundaries in panels
(b)--(c).  The variance uses $Q-1$ in the denominator.}
\label{fig:star-recon}
\end{figure}

Figure~\ref{fig:star-recon} shows one star realization, the ensemble mean, and
the unbiased variance.  The mean retains the common high-indicator
obstacle footprint, while the
variance is concentrated along the moving boundary and is strongest where the
prior produces the largest normal displacement.

\begin{figure}[tbhp]
\centering
\includegraphics[width=0.6\textwidth]{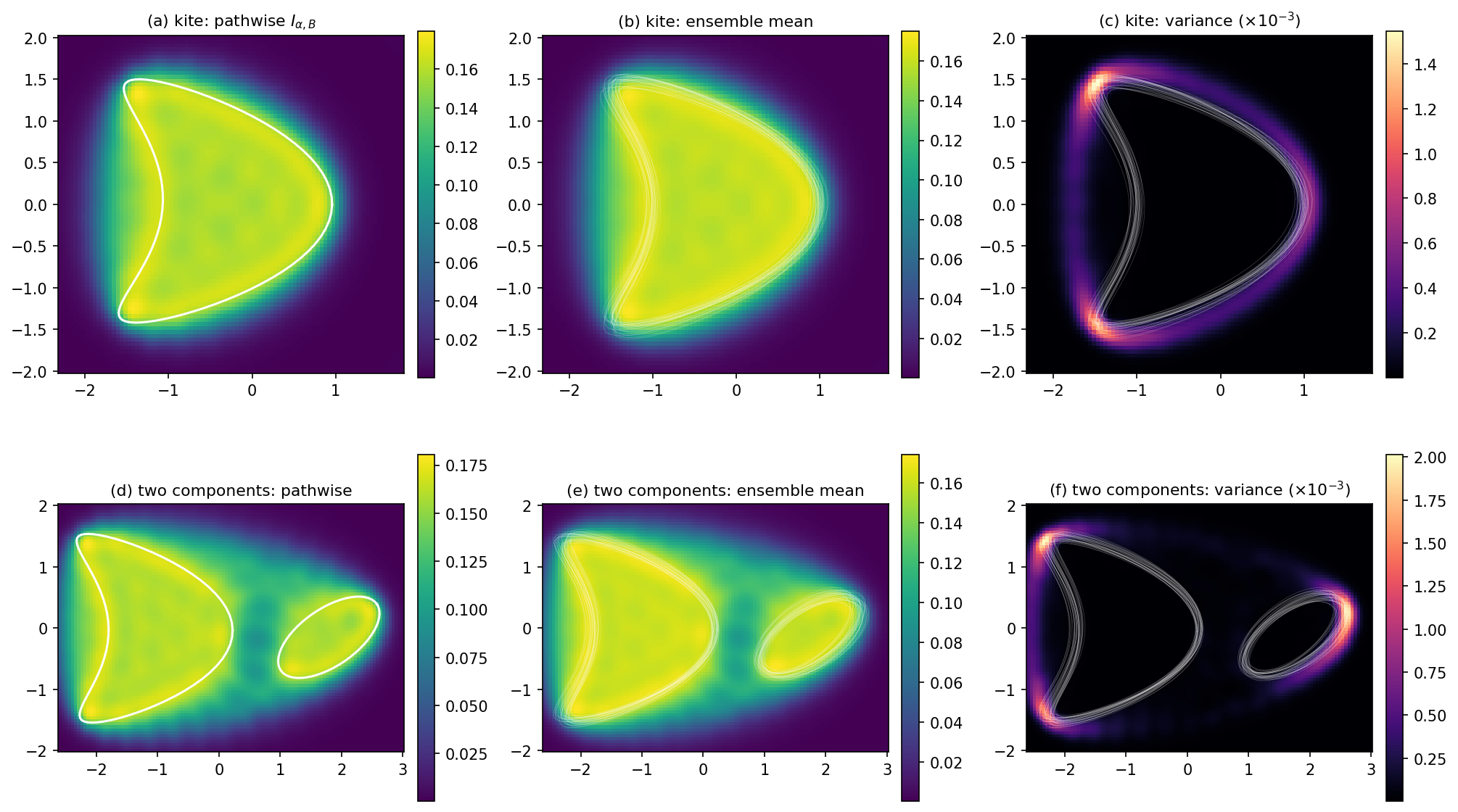}
\caption{Pathwise, mean, and unbiased-variance fields for the admissible affine
kite prior (top) and the two-component configuration (bottom).}
\label{fig:kite-recon}
\end{figure}

Figure~\ref{fig:kite-recon} gives the corresponding nonconvex and
multi-component examples.  The affine prior preserves the reentrant kite
geometry, simplicity, and component separation in every realization.

\subsection{Repeated measurement-noise study}
\begin{figure}[tbhp]
\centering
\includegraphics[width=0.6\textwidth]{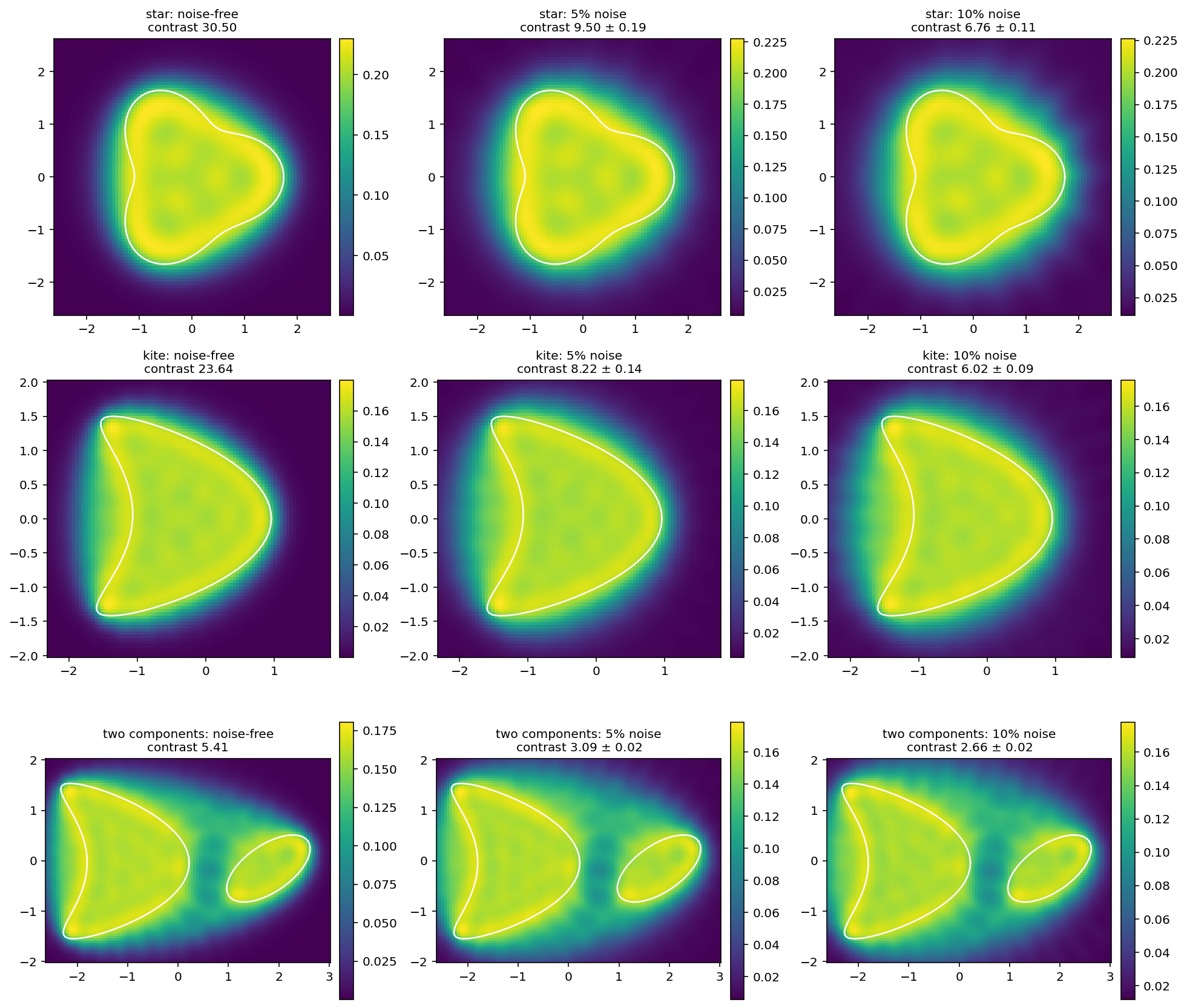}
\caption{Indicator fields under measurement noise.  Titles report the
$20$-draw mean and standard deviation of the contrast at nonzero noise
levels.}
\label{fig:noise}
\end{figure}

For a fixed representative realization of each geometry, let $G$ have
independent standard circular complex Gaussian entries and perturb the raw
sampled matrix $K_M$ of \eqref{eq:weighted-matrix}, before quadrature
weighting, by
\begin{equation}\label{eq:noise-level}
 E_\delta
 =\delta\,\norm{K_M}_{F}\frac{G}{\norm{G}_{F}},
 \qquad
 \delta\in\{0,0.05,0.10\};
\end{equation}
under uniform weights the weighted matrix and its perturbation are both
multiplied by $2\pi/M$, so the relative level is the same on either scale.
Thus the \emph{direction} is isotropic Gaussian and the relative Frobenius
norm is fixed exactly; after normalization the entries of $E_\delta$ are not
independent, so this experiment is not a realization of the
independent-frame model of Section~\ref{sec:discrete}.  Every nonzero level
uses $20$ independent matrix directions.  The contrast is the mean indicator
over the true interior divided by the mean over exterior grid points at
distance at least $0.55$ from the boundary; Table~\ref{tab:noise} reports
mean $\pm$ sample standard deviation and Figure~\ref{fig:noise} shows one
representative draw per level.

\begin{table}[tbhp]
\tbl{Interior/exterior indicator contrast over repeated noise draws.}
{\begin{tabular}{lrrr}
\toprule
geometry & noise-free & 5\% noise & 10\% noise\\
\midrule
star & 30.50 & $9.50\pm0.19$ & $6.76\pm0.11$\\
kite & 23.64 & $8.22\pm0.14$ & $6.02\pm0.09$\\
two components & 5.41 & $3.09\pm0.02$ & $2.66\pm0.02$\\
\bottomrule
\end{tabular}}
\label{tab:noise}
\end{table}

The draw-to-draw spread is small ($1$--$2\%$ of the mean), so across these
twenty sampled noise directions the degradation was stable rather than an
artefact of a particular draw.  The two-component contrast is lower at every level
because the exterior reference region at distance $0.55$ intersects the
channel between the components, where the indicator is elevated, and the
$B$-clipping flattens the interior mean.

\subsection{Random-disc UQ and variance leakage}
\begin{figure}[tbhp]
\centering
\includegraphics[width=0.58\textwidth]{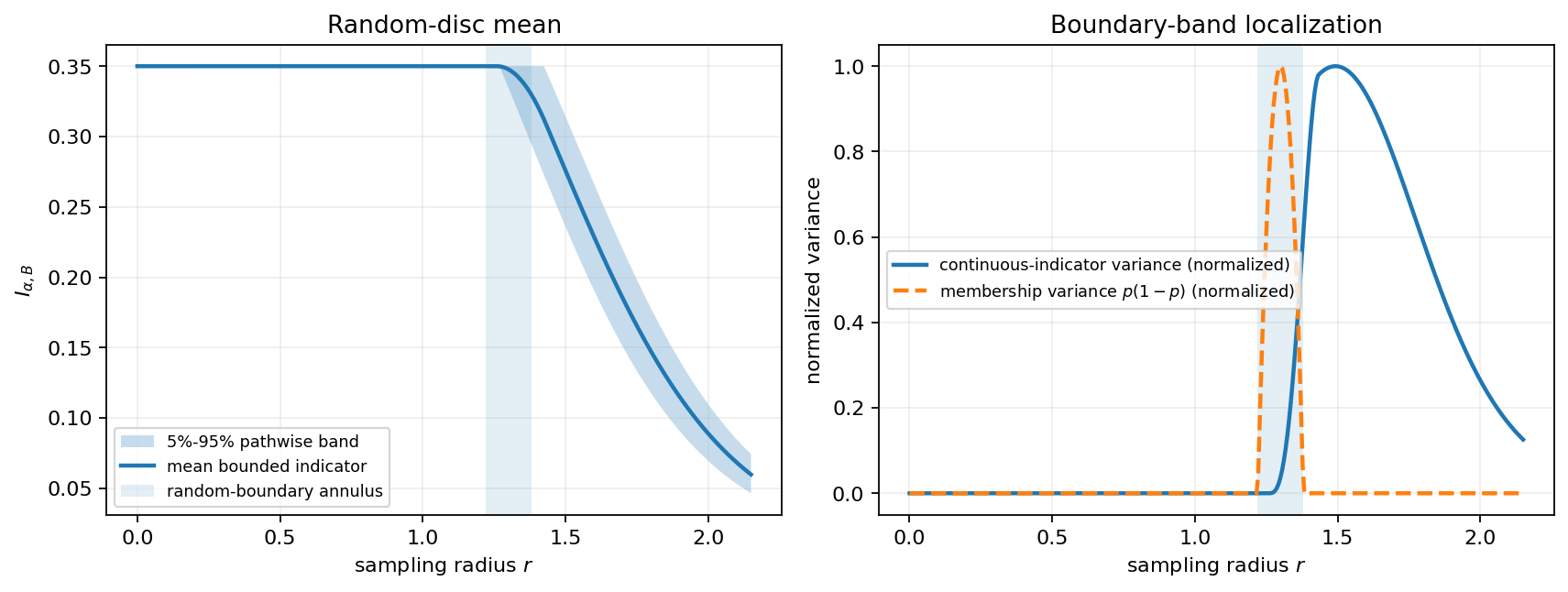}
\caption{Random-disc UQ.  The membership variance has exact boundary-band
support; the continuous regularized indicator has controlled but nonzero
leakage.}
\label{fig:disc-uq}
\end{figure}

\looseness=-1
Finally, let the disc radius be uniform on $[1.22,1.38]$, take
$\alpha=10^{-4}$, $B=0.35$, and evaluate the regularized Picard field of the
impedance disc from the eigenvalues \eqref{eq:disc-kappa}, namely
$W_{\alpha,L}(r)=2\pi\sum_{|m|\le L}J_m(kr)^2\gamma_m(\gamma_m+\alpha)^{-2}$
truncated at $L=24$.  Figure~\ref{fig:disc-uq} compares the continuous bounded
indicator with the exact membership field.  The Bernoulli membership variance
is supported exactly in the random-radius annulus.  The continuous-indicator
variance is concentrated near that annulus but is broadened and shifted by the
regularized imaging kernel.  This is consistent with the leakage mechanism bounded by
Proposition~\ref{prop:variance-leakage}; that proposition controls magnitude,
not the exact displacement or width of the variance peak.  An exact support
claim is therefore made only for membership or thresholded fields.

\subsection{Coverage and frame-averaging study}\label{sec:coverage-num}
\begin{figure}[tbhp]
\centering
\includegraphics[width=0.98\textwidth]{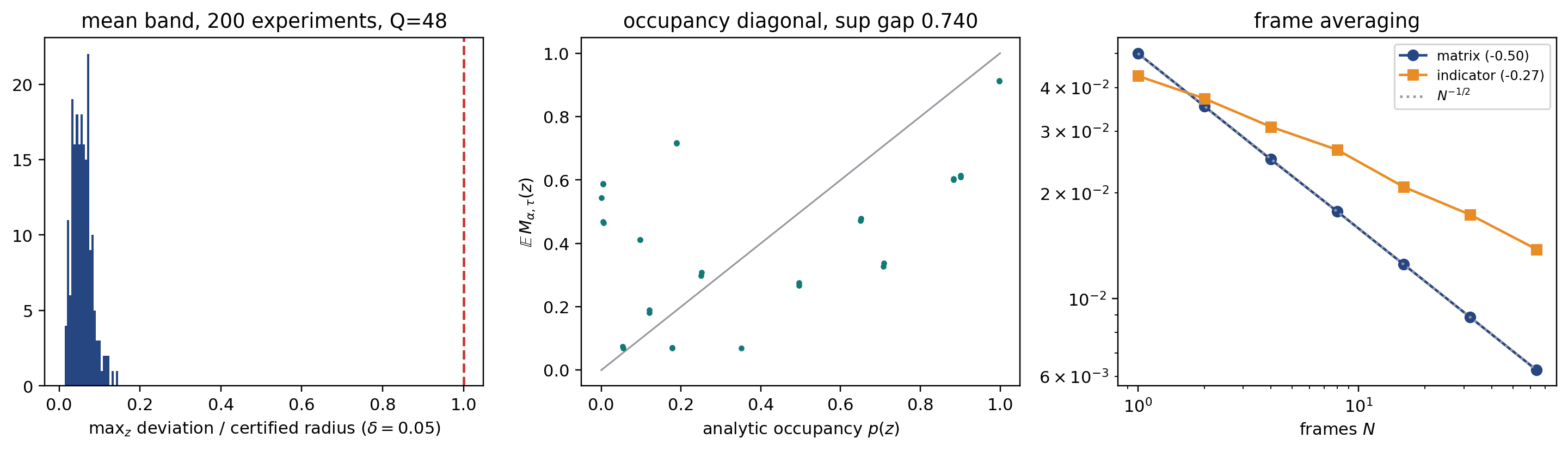}
\caption{Coverage study on the star prior.  Left: distribution over $200$
disjoint $Q=48$ experiments of the worst grid deviation of the empirical
mean, as a fraction of the $\delta=0.05$ radius of
Theorem~\ref{thm:grid-concentration}.  Middle:
calibrated occupancy frequency against the analytic prior occupancy on the
uncertain grid points.  Right: median frame-averaging error of the data
matrix and of the supremum indicator deviation.}
\label{fig:coverage}
\end{figure}

\looseness=-1
Throughout this subsection \emph{coverage} denotes the empirical frequency
with which a stated bound contains the quantity it is meant to bound, over
independent repetitions of the experiment.
The statistical guarantees of Section~\ref{sec:uq} are examined
empirically, against a high-accuracy Monte Carlo reference, on the star
prior with the combined-field solver, on a fixed $21\times21$
grid ($J=441$) over the sampling window, at the ensemble parameters $k=4$,
$M=64$, $\alpha=1.9635\times10^{-4}$, $B=0.35$.  A pool of $14400$
independent geometry draws is split into a reference part ($4800$ paths)
and $200$ disjoint experiments of $Q=48$ draws each; the indicators are
noise-free, matching the geometry-sampling setting of
Theorem~\ref{thm:grid-concentration}.  The reference paths provide
reference estimates of $\E X_1$ and $\Var X_1$, with pointwise standard
errors at most $B/(2\sqrt{4800})\approx2.5\times10^{-3}$ and
$B^{2}/(2\sqrt{4800})\approx8.8\times10^{-4}$, about $2\%$ of the
corresponding radii of Theorem~\ref{thm:grid-concentration} below.  All
deviations are measured against this reference, so the reported coverage is
an empirical illustration rather than an exact verification.  Figure~\ref{fig:coverage} summarizes
the three groups of results reported below.

\looseness=-1
\emph{Simultaneous bands.}  At $\delta=0.05$ every one of the $200$
experiments satisfies both the mean and the variance band; the largest
observed values of $\max_z|\overline X_Q-\E X_1|$ and
$\max_z|S_Q^2-\Var X_1|$ are $0.146$ and $0.015$ times the radii of
Theorem~\ref{thm:grid-concentration}.  Even at the nominal level $\delta=0.5$ the empirical coverage
remains $200/200$ (worst ratios $0.167$ and $0.017$).  The
distribution-free bounds are valid and, as expected, conservative ---
by roughly a factor $7$ for the mean and, at the small variance scale of
this ensemble ($\Var X_1\le10^{-3}$ against radius $0.039$), by much more
for the variance.

\looseness=-1
\emph{Occupancy.}  The threshold is calibrated on the reference part by
minimizing misclassification against the analytic star membership
($z\in D_Y\iff|z|<r_Y(\theta_z)$), giving $\tau^*=5.06$ and a held-out
misclassification rate of $2.49\%$.  The simultaneous Hoeffding band for
the occupancy frequency covers in $200/200$ experiments, and per-point
$95\%$ Hoeffding intervals on the $30$ uncertain grid points
($0.05<\E M<0.95$) cover in $99.6\%$ of cases; Wald intervals are slightly
anticonservative ($91.3\%$), as expected.  Against the prior-law occupancy
$p(z)$, computed from $2\times10^6$ direct draws of the radius law, the
fixed-$(\alpha,\tau)$ bias $|\E M-p|$ has mean $0.025$ over the grid but
supremum $0.74$, concentrated in the boundary annulus ($22$ of $441$
points exceed $0.2$).  Decreasing the regularization by two decades, with
recalibrated threshold, lowers the misclassification to $1.97\%$ and
shrinks the badly biased set (points above $0.5$ drop from $7$ to $2$;
mean gap $0.020$) while the worst annulus point barely moves.  This pointwise behavior is
consistent with the logarithmically slow diagonal convergence predicted by the severe smoothing of the far-field operator
\cite{PartII}, while the
measure-level error is already small.  At fixed $(\alpha,\tau)$ the
quantity controlled by the intervals is therefore $\E M_{\alpha,\tau}$, and its transfer to
$p$ is reliable away from the boundary band.

\looseness=-1
\emph{Nested components.}  Under an entrywise complex Gaussian frame model
(a sub-Gaussian instance of the frame hypothesis of
Section~\ref{sec:discrete}; single-frame relative Frobenius level $5\%$),
a nested reference ensemble of $1200$ geometries with $24$ replicates each
provides reference estimates of $V_{\rm meas}$ and $V_{\rm geom}$; $100$ experiments with $Q=12$
groups and $R=4$ replicates give component estimates whose averages differ
from the reference by at most $5\times10^{-7}$ and $6\times10^{-5}$, against
component scales $8\times10^{-6}$ and $7\times10^{-4}$.  These differences are
consistent with the exact unbiasedness of
Theorem~\ref{thm:nested-variance}.  The bounds of
Theorem~\ref{thm:nested-concentration} cover in $100/100$ experiments,
with worst deviation-to-radius ratios below $10^{-2}$.

\looseness=-1
\emph{Frame averaging.}  For $N$-frame averages, the median relative
operator error follows the $N^{-1/2}$ law with fitted slope $-0.500$ over
$N\in\{1,\dots,64\}$ (relative error $5.0\%$ at $N=1$ down to $0.63\%$ at
$N=64$).  The supremum indicator error at the fixed ensemble $\alpha$
decays more slowly (fitted tail slope $-0.31$; $0.043$ down to $0.014$)
than the relative data-matrix error.  This is
consistent with spectral amplification by the severely smoothing
far-field operator, and with the need to couple $\alpha$ to an absolute
weighted data-error bound as in \eqref{eq:joint-coupling} rather than hold
it fixed.

\subsection{Scope of the numerical evidence}
\looseness=-1
The experiments validate the weighted, phase-rotated implementation, the
analytic disc formulas, the layer-operator jump conventions --- including
the combined-field blocks $K_k$ and $N_k$ --- reciprocity, multiple
scattering, angular resolution, resonance-free ensemble generation with
recorded pathwise conditioning, and the qualitative behavior of
pathwise and ensemble indicators.  Section~\ref{sec:coverage-num} adds finite-grid coverage for the mean and
variance bands (Theorem~\ref{thm:grid-concentration}), the thresholded
occupancy construction, the nested variance components
(Theorem~\ref{thm:nested-concentration}), and the $N^{-1/2}$
frame-averaging law, all on the star prior.  Still untested are the
continuum band of Theorem~\ref{thm:continuum-band} beyond its
finite-grid core, and systematic sensitivity to $B$ and the prior beyond
the two regularization levels reported.

\FloatBarrier

\section{Conclusions and scope}\label{sec:conclusion}

\looseness=-1
The pathwise stochastic factorization developed here rests on a
normalization-correct data-to-pattern factorization.  In the physical
far-field convention of \eqref{eq:far-asymptotic},
$e^{-i\pi/4}F=-\sqrt{8\pi k}\,GT^*G^*$, so the positive scalar and the
two-dimensional phase must be retained in the proof even though they do not
change the final range.  Positive real impedance removes the direct and
adjoint real interior impedance eigenvalues.  The transferred optical
identity gives $\Imag(e^{-i\pi/4}F)\succeq0$ with strict quadratic-form
positivity.  The exact range operator is
$F_\#=|\Real(e^{-i\pi/4}F)|+\Imag(e^{-i\pi/4}F)$; the second absolute value
in the empirical matrix is a noise stabilization.

\looseness=-1
The finite-sample analysis matches the weighted Hilbert-space
discretization.  Angular approximation, frame noise, test-vector error, and
regularization are kept separate, deterministic confidence schedules specify
the probability mode, and, under uniform geometry-family hypotheses,
Theorem~\ref{thm:end-to-end} combines these errors with geometry sampling
into convergence of empirical means and variances.  The regularization rates that this analysis
calls for, and the sharp disc asymptotics that determine them, are
studied separately in \cite{PartII}.

\looseness=-1
The uncertainty analysis goes beyond a pointwise ensemble variance.  Bounded
indicators give finite-grid concentration and continuum bands for the mean
regularized indicator.  A diagonal threshold theorem recovers the physical
occupancy probability.  A nested design separates between-geometry from
within-geometry variability with unbiased estimators and concentration
bounds; the guarantees
admitting finite-grid tests are examined empirically, against a
high-accuracy Monte Carlo reference, in the coverage study of
Section~\ref{sec:coverage-num}.

\looseness=-1
Several limitations are explicit.  The single-layer representation can
encounter representation-induced interior Dirichlet resonances; the pathwise
recomputation in Table~\ref{tab:conditioning-ensemble} exhibits condition numbers up
to $1.99\times10^3$ on those seeds, and the worst case grows to
$1.1\times10^5$ over longer random runs while the combined field stays
below $6.8$ on the identical draws
(Table~\ref{tab:conditioning-tail}).  The reported ensembles and the
repeated-noise study are therefore generated with the validated
combined-field solver \eqref{eq:combined-field-bie}, which the same
recomputation shows to be well conditioned on every sampled path, the conditioning of
each realization being recorded alongside it; uniform discrete conditioning
over the whole family is not claimed, the collective-compactness question
being open.  Genuine limited aperture requires a separate restriction or
data-recovery analysis; more general complex or reactive impedances require
rechecking the sign and adjoint-eigenvalue hypotheses; and the
present UQ is a prior push-forward rather than a posterior conditioned on
one observed data set \cite{Stuart10}.

\looseness=-1
The coverage study empirically supports the finite-grid band, occupancy, and
variance-component guarantees and the $N^{-1/2}$ frame law on the star
prior; the continuum bands remain to be tested beyond their finite-grid
core, and on the solver side the regularized second-kind combined-field
variant, with the collective-compactness verification it requires, is the
natural next step.  Mathematically, the sharp spectral asymptotics
of \cite{PartII} suggest seeking geometry-uniform counterparts.

\section*{Disclosure statement}
No potential conflict of interest was reported by the authors.
\section*{Funding}
No funding was received.

\bibliographystyle{tfnlm}
\bibliography{references}

\appendix

\section{Solver validation}\label{app:validation}

This appendix and Appendix~\ref{app:params} collect the validation
of the Nystr\"om discretizations \cite{Kress14} used in the main text,
together with a complete specification of the parameters, geometries,
grids, seeds, and noise scaling used in every experiment.  The present
appendix reports agreement with the analytic disc far field and with the
Fourier eigenvalues of the layer operators, reciprocity, boundary-node and
two-component refinement, and the far-separation interaction defect; the
analytic disc quantities are those of \eqref{eq:disc-am} and
\eqref{eq:disc-layer-eigs}.

Table~\ref{tab:solver-validation} reports agreement near machine precision for the
single-disc far field and for the Fourier modes $m=0,1,2,3$ of the four
layer operators entering the two formulations.  These are the single layer $S_k$ and
adjoint double layer $K_k'$ of the single-layer equation, and the double
layer $K_k$ and hypersingular operator $N_k$ of the combined field.  The
parameters are $a=1.3$, $k=\lambda=1$, and $32$ directions.  The analytic
eigenvalues are $s_m$ and $\kappa_m'$ of \eqref{eq:disc-layer-eigs}
together with
\begin{equation}\label{eq:disc-KN-eigs}
 \kappa_m=-\frac12+\frac{i\pi ka}{2}J_m'(ka)H_m^{(1)}(ka),
 \qquad
 \nu_m=\frac{i\pi k^2a}{2}J_m'(ka)H_m^{(1)\prime}(ka),
\end{equation}
for $K_k$ and $N_k$ respectively, in the exterior jump convention of
\eqref{eq:bie}; $N_k$ is assembled through Maue's identity from the same
Kress-quadrature single layer \cite{Kress14}.

\begin{table}[tbhp]
\tbl{Disc and layer-operator validation.  Maximum error over the Fourier
modes $m=0,1,2,3$.}
{\small
\setlength{\tabcolsep}{4pt}
\begin{tabular}{rrrrrr}
\toprule
boundary nodes $n$ & far-field & $S_k$ & $K_k'$ & $K_k$ & $N_k$\\
\midrule
64  & $5.07\times10^{-16}$ & $3.74\times10^{-16}$ & $1.47\times10^{-15}$
    & $3.34\times10^{-16}$ & $8.66\times10^{-15}$\\
128 & $5.99\times10^{-16}$ & $3.19\times10^{-16}$ & $3.22\times10^{-15}$
    & $1.21\times10^{-15}$ & $1.98\times10^{-14}$\\
256 & $6.06\times10^{-16}$ & $3.51\times10^{-16}$ & $1.09\times10^{-15}$
    & $1.19\times10^{-15}$ & $6.26\times10^{-14}$\\
\bottomrule
\end{tabular}}
\label{tab:solver-validation}
\end{table}

Reciprocity is checked through
the relation $u^\infty(\widehat x,d)=u^\infty(-d,-\widehat x)$ \cite{CK19}.
Relative errors are
$1.87\times10^{-15}$ for the kite,
$1.20\times10^{-15}$ for a tilted ellipse, and
$2.48\times10^{-13}$ for the two-component configuration; the larger
multi-component value reflects the plain trapezoidal quadrature of the smooth
cross-interaction blocks at finite separation.  The optical identity
\eqref{eq:optical} is verified mode by mode on the disc
($a=1.3$, $k=\lambda=1$), where the largest residual over the first twenty-five
Fourier modes is $4.4\times10^{-16}$.

A boundary-node refinement study at $k=7$ compares the kite far field with a
$192$-node reference.

\begin{table}[tbhp]
\tbl{Boundary-node convergence for the kite at $k=7$ and $M=96$.}
{\small
\begin{tabular}{rrrrrr}
\toprule
$n$ & 48 & 64 & 80 & 96 & 128\\
\midrule
relative error
& $6.85\times10^{-2}$
& $3.15\times10^{-4}$
& $1.21\times10^{-7}$
& $6.02\times10^{-10}$
& $1.75\times10^{-12}$\\
\bottomrule
\end{tabular}}
\label{tab:boundary-convergence}
\end{table}
Thus the $96$-node kite discretization used below is well converged at the
selected frequency.  For the two-component scene, a
combined-field refinement against a $(320,224)$-node reference gives
relative errors $2.07\times10^{-8}$, $6.79\times10^{-12}$,
$5.31\times10^{-15}$, and $1.14\times10^{-15}$ at the node pairs
$(80,56)$, $(120,84)$, $(160,112)$, and $(240,168)$.  The relative
far-field change at the node counts used in the ensembles is
$2.07\times10^{-8}$, so the reported reconstructions are numerically
converged with respect to boundary refinement.

For two separated components, translation does not make the isolated far field
of the remote component small in operator norm.  Writing $K_{1\cup2}$,
$K_1$, $K_2$ for the sampled far-field matrices of the pair and of each
component in isolation, the appropriate interaction defect is
\begin{equation}\label{eq:interaction-defect}
 \frac{\norm{K_{1\cup2}-K_1-K_2}_{F}}{\norm{K_1}_{F}},
\end{equation}
where $d$ denotes the center-to-center separation of the two components.
Table~\ref{tab:separation} approaches the expected two-dimensional
$d^{-1/2}$ multiple-scattering decay.

\begin{table}[tbhp]
\tbl{Far-separation interaction defect.}
{\begin{tabular}{rrrrrr}
\toprule
separation $d$ & 5 & 10 & 20 & 40 & 80\\
\midrule
interaction defect
& 0.211531 & 0.154911 & 0.111814 & 0.079910 & 0.056825\\
local decay exponent
& -- & 0.449 & 0.470 & 0.485 & 0.492\\
\bottomrule
\end{tabular}}
\label{tab:separation}
\end{table}

\section{Complete numerical specification}\label{app:params}

All quantities needed to reproduce the experiments of the main text are
collected here.  Throughout, the impedance is the constant $\lambda=1$; the
combined-field coupling is $\eta=k$ in every ensemble; boundary curves are
sampled at equispaced parameter values; and all pseudorandom draws use
\texttt{numpy} PCG64 generators seeded as listed.

\paragraph*{Nominal geometries.}
With $t_j=2\pi j/n$, the reference curves are
\begin{align*}
 \text{star:}&\quad r(t)=1.45+0.34\cos 3t
   +0.10\!\sum_{m=1}^{4}\frac{\xi_m\cos mt+\eta_m\sin mt}{m},\\
 \text{kite:}&\quad x(t)=s\bigl(\cos t+0.65\cos 2t-0.65,\ 1.5\sin t\bigr)+c,\\
 \text{ellipse:}&\quad x(t)=R_\vartheta\,(a_1\cos t,\ a_2\sin t)^{\!\top}+c,
\end{align*}
with $\xi_m,\eta_m\sim\mathcal U[-1,1]$ and $R_\vartheta$ the rotation
by $\vartheta$.  The kite uses $s=1$, $c=(0,0)$
in the single-obstacle runs.  The two-component scene uses the kite with
$c=(-0.8,0)$ and the ellipse with $a_1=0.95$, $a_2=0.45$,
$\vartheta=0.6$, $c=(1.75,-0.10)$.

\paragraph*{Ensembles.}

\begin{center}
\small
\setlength{\tabcolsep}{3pt}
\begin{tabular}{lccccccl}
\toprule
scene & $k$ & $M$ & boundary nodes & $Q$ & $\alpha$ & $B$ & sampling grid\\
\midrule
star & 4 & 64 & 96 & 48 & $1.9635\times10^{-4}$ & 0.35
 & $90^2$ on $[-2.6,2.6]^2$\\
kite & 7 & 96 & 96 & 48 & $7.8540\times10^{-5}$ & 0.35
 & $90\times86$ on $[-2.3,1.8]\times[-2,2]$\\
two comp. & 7 & 96 & $80+56$ & 36 & $7.8540\times10^{-5}$ & 0.35
 & $110\times84$ on $[-2.6,3]\times[-2,2]$\\
\bottomrule
\end{tabular}
\end{center}

Geometry seeds are consecutive integers beginning at $3000$ (star),
$11\,000$ (kite) and $22\,000$ (two components); the pathwise panels of the
reconstruction figures show the first seed of each block.

\paragraph*{Repeated-noise study.}
Applied to the first realization of each scene, with $E_\delta$ as in
\eqref{eq:noise-level} and $20$ draws per nonzero level.  The draw seeds are
$700\,000+10^4 i+100\lfloor100\delta\rfloor+d$, with $i=0,1,2$ indexing the
scene and $d=0,\dots,19$ the draw.  The exterior reference region excludes
points within $0.55$ of the boundary.

\paragraph*{Random-disc study.}
Radius uniform on $[1.22,1.38]$, $k=\lambda=1$, $\alpha=10^{-4}$, Fourier
cutoff $|m|\le24$, $B=0.35$; the indicator is evaluated from the analytic
eigenvalues \eqref{eq:disc-am} and \eqref{eq:disc-kappa} rather than from a
discretized solve.

\paragraph*{Coverage study.}
Star prior with $k=4$, $M=64$, $\alpha=1.9635\times10^{-4}$, $B=0.35$,
$96$ boundary nodes, on a $21\times21$ grid ($J=441$) over
$[-2.6,2.6]^2$.  The path pool is $14\,400$ draws with seeds from
$500\,000$, of which the first $4800$ form the reference part and the remainder $200$
disjoint experiments of $Q=48$.  The nested design uses $1200$ reference
geometries (seeds from $600\,000$, $24$ replicates each) and $100$
experiments of $Q=12$ groups with $R=4$ replicates (seeds from
$7\,000\,000$); replicate noise seeds are offset by $10^7$ and $2\times10^7$
respectively.  The frame-averaging sweep uses the realization with seed
$800\,001$ and $40$ draws per frame count $N\in\{1,2,4,8,16,32,64\}$
(seeds from $800\,100$).

For a complex matrix $K\in\mathbb C^{M\times M}$, one noisy frame is
$K+\sigma G$ with $G$ having independent standard complex Gaussian entries
(unit variance per real component) and
\begin{equation}\label{eq:frame-sigma}
 \sigma=\frac{\varrho\,\norm{K}_F}{\sqrt2\,M},
 \qquad \varrho=0.05 ,
\end{equation}
so that $\E\norm{\sigma G}_F=\varrho\norm{K}_F$ up to $O(M^{-1})$.  The
single-frame relative Frobenius noise level is therefore $5\%$.  An $N$-frame average
draws $N$ independent such matrices and averages them.

\end{document}